\input amstex\documentstyle{amsppt}  
\pagewidth{12.5 cm}\pageheight{19 cm}\magnification\magstep1
\topmatter
\title Unipotent elements in small characteristic\endtitle
\author G. Lusztig\endauthor
\address Department of Mathematics, M.I.T., Cambridge, MA 02139\endaddress
\dedicatory{Dedicated to Vladimir Drinfeld on the occasion of his 50-th birthday}
\enddedicatory
\thanks Supported in part by the National Science Foundation.\endthanks
\endtopmatter   
\document
\define\Nil{\text{\rm Nil}}
\define\lar{\la,\ra}

\define\part{\partial}

\define\dx{\dot x}
\define\dy{\dot y}

\define\mpb{\medpagebreak}

\define\da{\dagger}

\define\Lie{\text{\rm Lie }}

\define\frl{\forall}
\define\pe{\perp}
\define\si{\sim}

\define\sqc{\sqcup}

\define\qua{\quad}

\define\bg{\bar g}

\define\bC{\bar C}
\define\bE{\bar E}

\define\bG{\bar G}

\define\bN{\bar N}

\define\bS{\bar S}
\define\bT{\bar T}
\define\bX{\bar X}
\define\bY{\bar Y}
\define\bZ{\bar Z}
\define\bV{\bar V}

\define\bbq{\bar{\QQ}_l}
\define\baP{\bar\Pi}

\define\op{\oplus}
   
\redefine\sp{\spadesuit}
\define\em{\emptyset}
\define\imp{\implies}
\define\ra{\rangle}
\define\n{\notin}
\define\iy{\infty}
\define\m{\mapsto}
\define\do{\dots}
\define\la{\langle}

\define\lra{\leftrightarrow}

\define\sm{\smallmatrix}
\define\esm{\endsmallmatrix}
\define\sub{\subset}

\define\T{\times}
\define\ti{\tilde}
\define\nl{\newline}
\redefine\i{^{-1}}

\define\ov{\overline}
\define\ot{\otimes}

\define\bcc{\thickfracwithdelims[]\thickness0}

\define\Ad{\text{\rm Ad}}
\define\Hom{\text{\rm Hom}}
\define\End{\text{\rm End}}

\define\Ind{\text{\rm Ind}}

\redefine\Im{\text{\rm Im}}

\define\a{\alpha}
\redefine\b{\beta}

\define\g{\gamma}
\redefine\d{\delta}
\define\e{\epsilon}

\redefine\o{\omega}
\define\p{\pi}
\define\ph{\phi}

\define\r{\rho}
\define\s{\sigma}
\redefine\t{\tau}
\define\th{\theta}

\redefine\l{\lambda}
\define\z{\zeta}
\define\x{\xi}

\redefine\G{\Gamma}
\redefine\D{\Delta}

\define\Si{\Sigma}
\define\Th{\Theta}

\define\Ph{\Phi}

\define\kk{\bold k}

\define\nn{\bold n}

\define\CC{\bold C}

\define\FF{\bold F}

\define\NN{\bold N}

\define\QQ{\bold Q}

\define\VV{\bold V}
\define\WW{\bold W}
\define\ZZ{\bold Z}

\define\ca{\Cal A}
\define\cb{\Cal B}
\define\cc{\Cal C}

\define\ce{\Cal E}
\define\cf{\Cal F}
\define\cg{\Cal G}

\define\cj{\Cal J}

\define\cl{\Cal L}
\define\cm{\Cal M}

\define\co{\Cal O}

\define\cs{\Cal S}

\define\cu{\Cal U}
\define\cv{\Cal V}

\define\cz{\Cal Z}
\define\cx{\Cal X}

\define\fF{\frak F}

\define\fP{\frak P}
\define\fQ{\frak Q}

\define\fS{\frak S}

\define\ta{\ti a}
\define\tb{\ti b}
\define\tc{\ti c}

\define\tih{\ti h}
\define\tj{\ti j}

\define\tu{\ti u}

\define\tB{\ti B}

\define\tD{\ti D}
\define\tE{\ti E}

\define\tG{\ti G}
\define\tH{\ti H}

\define\tN{\ti N}

\define\tU{\ti U}

\define\tW{\ti W}

\define\tZ{\ti Z}

\define\tix{\ti\x}

\define\Ir{\text{\rm Irr}}
\define\sps{\supset}

\define\che{\check}

\define\prim{\text{\rm prim}}    
\define\Bil{\text{\rm Bil}}
\define\bcb{\bar{\cb}}
\redefine\bcc{\bar{\cc}}
\define\gr{\text{\rm gr}}
\define\Symp{\text{\rm Symp}}
\define\hag{\hat\g}
\define\bla{\blacktriangle}
\redefine\l{\vartriangle}
\define\DLP{DLP}
\define\DE{D1}
\define\DEII{D2}
\define\EN{E}
\define\KA{Ka}
\define\KO{Ko}
\define\LC{L1}
\define\LN{L2}
\define\LS{LS}
\define\MI{M}
\define\SH{Sh}
\define\SPA{S1}
\define\SPAII{S2}
\define\SP{Sp}
\define\WA{W}
\head Introduction\endhead
\subhead 0.1\endsubhead
Let $\kk$ be an algebraically closed field of characteristic exponent $p\ge1$. Let $G$
be a reductive connected algebraic group over $\kk$. Let $\cu$ be the variety of 
unipotent elements of $G$. The unipotent classes of $G$ are the orbits of the
conjugation action of $G$ on $\cu$. The theory of Dynkin and Kostant \cite{\KO} 
provides a classification of unipotent classes of $G$ assuming that $p=1$. It is known 
that this classification remains valid when $p\ge2$ is assumed to be a good prime for 
$G$. But the analogous classification problem in the case where $p$ is a bad prime for 
$G$ is more complicated. In every case a classification of unipotent classes is known: 
see \cite{\WA} for classical groups and \cite{\EN,\SH,\MI} for exceptional groups; but 
from these works it is difficult to see the general features of the classification. 

One of the aims of this paper is to present a picture of the unipotent elements which
should apply for arbitrary $p$ and is as close as possible to the picture for $p=1$. 

In 1.4 we observe that the set of unipotent classes in $G$ can be parametrized by a set
$\cs^p(\WW)$ of irreducible representations of the Weyl group $\WW$ which can be 
described apriori purely in terms of the root system. This explains clearly why the 
classification is different for small $p$. 

In 1.1 we restate in a more precise form an observation of \cite{\LN} according to 
which $\cu$ is naturally partitioned into finitely many "unipotent pieces" which are 
locally closed subvarieties stable under conjugation by $G$; the classification of 
unipotent pieces is independent of $p$. For $p=1$ or a good prime, each unipotent piece
is a single conjugacy class. When $p$ is a bad prime a unipotent piece is in general a 
union of several conjugacy classes. Also each unipotent piece has 
some topological properties which are independent of $p$ (for example, over a finite 
field, the number of points of a unipotent piece is given by a formula independent of 
the characteristic). 

Another aim of this paper is the study of $\cb_u$, the variety of Borel subgroups of 
$G$ containing a unipotent element $u$. It is known \cite{\SP} that when $p$ is a good
prime, the $l$-adic cohomology spaces of $\cb_u$ are pure. We would like to prove a 
similar result in the case where $p$ is a bad prime. We present a method by which this 
can be achieved in a number of cases. Our strategy is to extend a technique from 
\cite{\DLP} in which (assuming that $p=1$), $\cb_u$ is analyzed by first partitioning 
it into finitely many smooth locally closed subvarieties using relative position of a 
point in $\cb_u$ with a canonical parabolic attached to $u$. Much of our effort is 
concerned with trying to eliminate reference to the linearization procedure of
Bass-Haboush (available only for $p=1$) which was used in an essential way in 
\cite{\DLP}. Our approach is based on a list of properties $\fP_1-\fP_8$ of unipotent 
elements of which the first five (resp. last three) are expected to hold in general 
(resp. in many cases). All these properties are verified for general linear and 
symplectic groups (any $p$) in \S2, \S3. In writing \S3 (on symplectic groups mostly 
with $p=2$) I found that the treatment in \cite{\WA} is not sufficient for this paper's
purposes; I therefore included a treatment which does not rely on \cite{\WA}.

{\it Notation.} When $p>1$ we denote by $\kk_p$ an algebraic closure of the field with 
$p$ elements. Let $\cb$ the variety of Borel subgroups of $G$. If $\G'$ is a subgroup 
of a group $\G$ and $x\in\G$ let $Z_{G'}(x)=\{z\in\G';zx=xz\}$. For a finite set $Z$ 
let $|Z|$ be the cardinal of $Z$. Let $l$ be a prime number invertible in $\kk$. For 
$a,b\in\ZZ$ let $[a,b]=\{z\in\ZZ;a\le z\le b\}$.

\head Contents\endhead
1. Some properties of unipotent elements.

2. General linear groups.

3. Symplectic groups.

4. The group $A^1(u)$.

5. Study of the varieties $\cb_u$.

\head 1. Some properties of unipotent elements\endhead
\subhead 1.1\endsubhead
$G$ acts naturally, by conjugation, on $\Hom(\kk^*,G)$ (homomorphisms of algebraic 
groups). The set of orbits $\Hom(\kk^*,G)/G$ is naturally in bijection with the 
analogous set $\Hom(\CC^*,G')/G'$ where $G'$ is a connected reductive group over $\CC$ 
of the same type as $G$. (Both sets may be identified with the set of Weyl group orbits
on the group of $1$-parameter subgroups of some maximal torus.) Let $\tD_{G'}$ be the 
set of all $\o\in\Hom(\CC^*,G')$ such that there exists a homomorphism of algebraic 
groups $\ti\o:SL_2(\CC)@>>>G'$ with $\ti\o\left(\sm t&0\\0&t\i\esm\right)=\o(t)$ for 
all $t\in\CC^*$. Now $\tD_{G'}$ is $G'$-stable; it has been described explicitly by 
Dynkin. Let $\tD_G$ be the unique $G$-stable subset of $\Hom(\kk^*,G)$ whose image in 
$\Hom(\kk^*,G)/G$ corresponds under the bijection 
$\Hom(\kk^*,G)/G\lra\Hom(\CC^*,G')/G'$ (as above) to the image of $\tD_{G'}$ in 
$\Hom(\CC^*,G')/G'$. Let $D_G$ be the set of sequences 
$\l=(G^\l_0\sps G^\l_1\sps G^\l_2\sps\do)$ of closed connected subgroups of $G$ such 
that for some $\o\in\tD_G$ we have (for $n\ge0$):

$\Lie G^\l_n=\{x\in\Lie G;\lim_{t\in\kk^*;t\to 0}t^{1-n}\Ad\o(t)x=0\}$.
\nl
Now $G$ acts on $D_G$ by conjugation and the obvious map $\tD_G@>>>D_G$ induces a 
bijection $\tD_G/G@>\si>>D_G/G$ on the set of orbits. If $\l\in D_G$ and $g\in G$ then 
$G^{g\l g\i}_n=gG^\l_ng\i$ for $n\ge0$. Moreover, $G^\l_0$ is a parabolic subgroup of 
$G$ with unipotent radical $G^\l_1$ and $G^\l_n$ is normalized by $G^\l_0$ for any $n$.
Moreover,

(a) $G^\l_2/G^\l_3$ is a commutative connected unipotent group;

(b) the conjugation action of $G^\l_0$ on $G^\l_2/G^\l_3$ factors through an action of 
$\bG^\l_0:=G^\l_0/G^\l_1$ on $G^\l_2/G^\l_3$.
\nl
Note also that $G^\l_n$ for $n\ne0,2$ are uniquely determined by $G^\l_0,G^\l_2$. 

Let $\bla$ be a $G$-orbit in $D_G$. Then $\tH^\bla:=\cup_{\l\in\bla}G^\l_2$ is a closed
irreducible subset of $\cu$ (since for $\l\in\bla$, $G^\l_2$ is a closed irreducible 
subset of $\cu$ stable under conjugation by $G^\l_0$ and $G/G^\l_0$ is projective). Let
$$H^\bla=\tH^\bla-\cup_{\bla'\in D_G/G;\tH^{\bla'}\subsetneqq\tH^\bla}\tH^{\bla'}.$$
For $\l\in D_G$ let $X^\l=G^\l_2\cap H^\bla$ where $\bla$ is the $G$-orbit of $\l$.    
Then $H^\bla$ is an open dense subset of $\tH^\bla$ stable under conjugation by $G$ and
$X^\l$ is an open dense subset of $G^\l_2$ stable under conjugation by $G^\l_0$. (We 
use that $D_G/G$ is finite.) Hence $H^\bla$ is locally closed in $\cu$. The subsets 
$H^\bla(\bla\in D_G/G)$ are called the {\it unipotent pieces} of $G$.

We state the following properties $\fP_1-\fP_5$.

$\fP_1$. {\it The sets $X^\l(\l\in D_G)$ form a partition of $\cu$.}

$\fP_2$. {\it Let $\bla\in D_G/G$. The sets $X^\l(\l\in\bla)$ form a partition of 
$H^\bla$. More precisely, $H^\bla$ is a fibration over $\bla$ with smooth fibres
isomorphic to $X^\l$ ($\l\in\bla$); in particular, $H^\bla$ is smooth.}

$\fP_3$. {\it The locally closed subets $H^\bla(\bla\in D_G/G)$ form a (finite) 
partition of $\cu$.}

$\fP_4$. {\it Let $\l\in D_G$. We have $G^\l_3X^\l=X^\l G^\l_3=X^\l$.}

$\fP_5$. {\it Assume that $\kk=\kk_p$. Let $F:G@>>>G$ be the Frobenius map 
corresponding to a split $\FF_q$-rational structure with $q-1$ sufficiently divisible. 
Let $\l\in D_G$ be such that $F(G^\l_n)=G^\l_n$ for all $n\ge0$ and let $\bla$ be the
$G$-orbit of $\l$. Then $|H^\bla(\FF_q)|,|X^\l(\FF_q)|$ are polynomials in $q$ with 
integer coefficients independent of $p$.}

Assume first that $p=1$ or $p\gg0$. By the theory of Dynkin-Kostant, for $\l\in D_G$ 
there is a unique open $G^\l_0$-orbit $X'{}^\l$ in $G^\l_2$; we then have a bijection 
of $D_G/G$ with the set of unipotent classes on $G$ which to the $G$-orbit $\bla$ of 
$\l\in D_G$ associates the unique unipotent class $H'{}^\bla$ of $G$ that contains 
$X'{}^\l$; moreover, if $g\in X'{}^\l$ then $Z_{G^\l_0}(g)=Z_G(g)$. As stated by
Kawanaka \cite{\KA}, the same holds when $p$ is a good prime of $G$ (but his argument 
is rather sketchy). To show that $\fP_1-\fP_3$ holds when $p$ is a good prime it then 
suffices to show that $X^\l=X'{}^\l$ for any $\l$. It also suffices to show that 
$X'{}^\l=G^\l_2\cap H'{}^\bla$ for $\l\in\bla$ as above. (Assume that 
$g\in G^\l_2\cap H'{}^\bla,g\n X'{}^\l$. Let $g'\in X'{}^\l$. By the definition of
$X'{}^\l$ and the irreducibility of $G^\l_2$, the dimension of the $G^\l_0$-orbit of 
$g$ is strictly smaller than the dimension of the $G^\l_0$-orbit of $g'$. Hence 
$\dim Z_{G^\l_0}(g)>\dim Z_{G^\l_0}(g')$. We have 
$\dim Z_G(g)\ge\dim Z_{G^\l_0}(g)$, $\dim Z_{G^\l_0}(g')=\dim Z_G(g')$ hence 
$\dim Z_G(g)>\dim Z_G(g')$. This contradicts the fact that $g,g'$ are $G$-conjugate.)
In this case we have $H^\bla=H'{}^\bla$ and $\tH^\bla$ is the closure of $H'{}^\bla$.

We expect that $\fP_1-\fP_5$ hold in general. In the case where $G=GL_n(\kk)$ (any $p$)
the validity of $\fP_1-\fP_5$ follows from 2.9. In the case where $G$ is a symplectic 
group (any $p$) the validity of $\fP_1-\fP_5$ follows from 3.13, 3.14. If $G$ is of 
type $E_n$ (any $p$) then one can deduce $\fP_1-\fP_5$ from the various lemmas in 
\cite{\MI}, or rather from the extensive computations (largely omitted) on which those 
lemmas are based; it would therefore be desirable to have an independent verification 
of these properties. The case of special orthogonal groups will be considered 
elsewhere. 

\mpb

We note the following consequence of $\fP_1$.

(c) {\it If $\l\in D_G$ and $u\in X^\l$ then $Z_G(u)\sub G^\l_0$.}
\nl
Let $g\in G$. Then $gug\i\in X^{g\l}$. Hence if $g\in Z_G(u)$ we have $u\in X^{g\l}$. 
Thus, $X^{g\l}\cap X^\l\ne\em$. From $\fP_1$ we see that $g\!\l=\l$. In particular 
$gG^\l_0g\i=G^\l_0$ and $g\in G^\l_0$, as required.

\subhead 1.2\endsubhead
Let $\l\in D_G$. We assume that $\fP_1-\fP_4$ hold for $\l$. Let 
$\p^\l:G^\l_2@>>>G^\l_2/G^\l_3$ be the obvious homomorphism. By $\fP_4$ we have 
$X^\l=(\p^\l)\i(\bX^\l)$ where $\bX^\l$ is a well defined open dense subset of 
$G^\l_2/G^\l_3$ stable under the action of $\bG^\l_0$. We wish to consider some 
properties of the sets $\bX^\l$ which may or may not hold for $G$. 

$\fP_6$. {\it If $u\in X^\l$ then $uG^\l_3=G^\l_3u$ is contained in the 
$G^\l_0$-conjugacy class of $u$. Hence $\g\m(\p^\l)\i(\g)$ is a bijection between the 
set of $\bG^\l_0$-orbits in $\bX^\l$ and the set of $G^\l_0$-conjugacy classes in 
$X^\l$.}

$\fP_7$. {\it Let $\g$ be a $\bG^\l_0$-orbit in $\bX^\l$. Let $\hag$ be the union of 
all $\bG^\l_0$-orbits in $\bX^\l$ whose closure contains $\g$. Thus, $\hag$ is an open 
subset of $\bX^\l$ and $\g$ is a closed subset of $\hag$. There exists a variety $\g_1$
and a morphism $\r:\hag@>>>\g_1$ such that the restriction of $\r$ to $\g$ is a finite
bijective morphism $\s:\g@>>>\g_1$ and the map of sets $\s\i\r:\hag@>>>\g$ is 
compatible with the actions of $\bG^\l_0$.}

$\fP_8$. {\it There exists a finite set $I$ and a bijection $J\m\Ph_J$ between the set 
of subsets of $I$ and the set of $G^\l_0$-orbits in $X^\l$ such that for any $J\sub I$,
the closure of $\Ph_J$ in $X^\l$ is $\cup_{J';J\sub J'}\Ph_{J'}$. Moreover, if $\kk,q$ 
are as in $\fP_5$ then there exists a function $I@>>>\{2,4,6,\do\},i\m c_i$ such that 
$|\Ph_J(\FF_q)|=\prod_{i\in J}(q^{c_i}-1)|\Ph_\em(\FF_q)|$ for any $J\sub I$.}

When $p=1$ or $p\gg0$ property $\fP_6$ can be deduced from the theory of 
Dynkin-Kostant; properties $\fP_7,\fP_8$ are trivial. In the case where $G=GL_n(\kk)$ 
(any $p$) the validity of $\fP_6$ follows from 2.9; properties $\fP_7,\fP_8$ are 
trivial. In the case where $G$ is a symplectic group (any $p$) the validity of 
$\fP_6-\fP_8$ follows from 3.14. $\fP_6$ is false for $G$ of type $G_2$, $p=3$. 

\subhead 1.3\endsubhead
Let $\VV$ be a finite dimensional $\QQ$-vector space. Let $R\sub\VV^*=\Hom(\VV,\QQ)$ be
a (reduced) root system, let $\che R\sub\VV$ be the corresponding set of coroots and 
let $\WW\sub GL(\VV)$ be the Weyl group of $R$. Let $\b\lra\che\b$ be the canonical 
bijection $R\lra\che R$. Let $\Pi$ be a set of simple roots for $R$ and let 
$\che\Pi=\{\che\a;\a\in\Pi\}$. Let 

$\Th=\{\b\in R;\b-\a\n R\qua\frl\a\in\Pi\}$,
$\ti\Th=\{\b\in R;\che\b-\che\a\n\che R\qua\frl\a\in\Pi\}$,

$\ti\ca=\{J\sub\Pi\cup\ti\Th;J\text{ linearly independent in }\VV^*\}$.
\nl
For any prime number $r$ let $\ca_r$ be the set of all $J\sub\Pi\cup\Th$ such that $J$ 
is linearly independent in $\VV^*$ and the torsion subgroup of 
$\sum_{\a\in\Pi}\ZZ\a/\sum_{\b\in J}\ZZ\b$ has order $r^k$ for some $k\in\NN$.

For any $J\in\ca_r$ or $J\in\ti\ca$ let $\WW_J$ be the subgroup of $\WW$ generated by
the reflections with respects to roots in $J$. For $W'=\WW$ or $\WW_J$ let $\Ir(W')$ be
the set of (isomorphism classes) of irreducible representations of $W'$ over $\QQ$. For
$E\in\Ir(W')$ let $b_E$ be the smallest integer $\ge0$ such that $E$ appears with 
non-zero multiplicity in the $b_E$-th symmetric power of $\VV$ regarded as a 
$W'$-module; if this multiplicity is $1$ we say that $E$ is good. If $J$ is as above 
and $E\in\Ir(W')$ is good then there is a unique $\tE\in\Ir(\WW)$ such that $\tE$ 
appears in $\Ind_{\WW_J}^\WW E$ and $b_{\tE}=b_E$; moreover, $\tE$ is good. We set 
$\tE=j_{\WW_J}^\WW E$.

Let $\cs_\WW\sub\Ir(\WW)$ be the set of special representations of $\WW$ (see 
\cite{\LC}). Now any $E\in\cs_\WW$ is good. Following \cite{\LC}, let $\cs^1_\WW$ be 
the set of all $E\in\Ir(\WW)$ such that $E=j_{\WW_J}^\WW E_1$ for some $J\in\ti\ca$ and
some $E_1\in\cs_{\WW_J}$. (Note that $\WW_J$ is like $\WW$ with the same $\VV$ and with
$R$ replaced by the root system with $J$ as set of simple roots; hence $\cs_{\WW_J}$ is
defined.) Now any $E\in\cs^1_\WW$ is good.

For any prime number $r$ let $\cs^r_\WW$ be the set of all $E\in\Ir(\WW)$ such that 
$E=j_{\WW_K}^\WW E_1$ for some $K\in\ca_r$ and some $E_1\in\cs^1_{\WW_K}$. (Note that 
$\WW_K$ is like $\WW$ with the same $\VV$ and with $R$ replaced by the root system with
$K$ as set of simple roots; hence $\cs^1_{\WW_K}$ is defined.) Now any $E\in\cs^1_\WW$ 
is good.

We have $\cs^1(\WW)\sub\cs^r(\WW)$. We have $\cs^1(\WW)=\cs^r(\WW)$ if $r$ is a good 
prime for $\WW$ and also in the following cases $\WW$ of type $G_2,r=2$; $\WW$ of type 
$F_4$, $r=3$; $\WW$ of type $E_6$; $\WW$ of type $E_7$, $r=3$; $\WW$ of type $E_8$, 
$r=5$. If $\WW$ is of type $G_2$ and $r=3$ then $\cs^r(\WW)-\cs^1(\WW)$ consists of a 
single representation of dimension $1$ coming under $j_{\WW_J}^\WW$ from a $\WW_J$ of 
type $A_2$. If $\WW$ is of type $F_4$ and $r=2$ then $\cs^r(\WW)-\cs^1(\WW)$ consists 
of four representations of dimensions $9/4/4/2$ coming under $j_{\WW_J}^\WW$ from a 
$\WW_J$ of type $C_3A_1/C_3A_1/B_4/B_4$. If $\WW$ is of type $E_7$ and $r=2$ then 
$\cs^r(\WW)-\cs^1(\WW)$ consists of a single representation of dimensions $84$ coming 
under $j_{\WW_J}^\WW$ from a $\WW_J$ of type $D_6A_1$. If $\WW$ is of type $E_8$ and 
$r=2$ then $\cs^r(\WW)-\cs^1(\WW)$ consists of four representations of dimensions 
$1050/840/168/972$ coming under $j_{\WW_J}^\WW$ from a $\WW_J$ of type 
$E_7A_1/D_5A_3/D_8/E_7A_1$. If $\WW$ is of type $E_8$ and $r=3$ then 
$\cs^r(\WW)-\cs^1(\WW)$ consists of a single representation of dimensions $175$ coming 
under $j_{\WW_J}^\WW$ from a $\WW_J$ of type $E_6A_2$. 

\subhead 1.4\endsubhead
Let $\WW$ be the Weyl group of $G$. Let $u$ be a unipotent element in $G$. Springer's 
correspondence (generalized to arbitrary characteristic) associates to $u$ and the 
trivial representation of $Z_G(u)/Z_G(u)^0$ a representation $\r_u\in\Ir(\WW)$. 
Moreover $u\m\r_u$ defines an injective map from the set of unipotent classes in $G$ to
$\Ir(\WW)$. Let $\cx^p(\WW)$ be the image of this map ($p$ as in 0.1). We state:

(a) {\it If $p=1$ we have $\cx^1(\WW)=\cs^1(\WW)$.} (See \cite{\LC}).

(b) {\it If $p>1$ we have $\cx^p(\WW)=\cs^p(\WW)$.} 
\nl
The proof of (b) follows from the explicit description of the Springer correspondence 
for small $p$ given in \cite{\LS},\cite{\SPAII}.

\head 2. General linear groups\endhead
\subhead 2.1\endsubhead
Let $\bcc$ be the category whose objects are $\ZZ$-graded $\kk$-vector spaces 
$\bV=\op_{a\in\ZZ}\bV_a$ such that $\dim\bV<\iy$; the morphisms are linear maps 
respecting the grading. Let $\bV\in\bcc$. For $j\in\ZZ$ let 
$\End_j(\bV)=\{T\in\Hom(\bV,\bV);T(\bV_a)\sub\bV_{a+j}\qua\frl a\}$. Let 
$\End_2^0(\bV)$ be the set of all $\nu\in\End_2(\bV)$ that satisfy the {\it Lefschetz 
condition}: $\nu^n:\bV_{-n}@>>>\bV_n$ is an isomorphism for any $n\ge0$. Let 
$\nu\in\End_2^0(\bV)$. Define a graded subspace $P^\nu=\bV^{\prim}$ of $\bV$ by 
$P^\nu_a=\{x\in\bV_a;\nu^{1-a}x=0\}$ for $a\le0$, $P^\nu_a=0$ for $a>0$. A standard 
argument shows that $N^{(a-c)/2}:P^\nu_c@>>>\bV_a$ is injective if 
$c\in a+2\ZZ,c\le a\le-c$ and we have

(a) $\op_{c\in a+2\ZZ;c\le a\le-c}P^\nu_c@>\si>>\bV_a,
(z_c)\m\sum_{c\in a+2\ZZ;c\le a\le-c}N^{(a-c)/2}z_c$.
\nl
We show:

(b) {\it Let $j\in\NN$, $R\in\End_{j+2}(\bV)$. Then $R=T\nu-\nu T$ for some
$T\in\End_j(\bV)$.}
\nl
Let $c\le0$. Since $\nu^{1-c}:\bV_{j-c}@>>>\bV_{j+c+2}$ is surjective, the induced map 

$\Hom(P^\nu_{-c},\bV_{j-c})@>>>\Hom(P^\nu_{-c},\bV_{j+c+2})$
\nl
is surjective. Hence there exists $\t_c\in\Hom(P^\nu_{-c},\bV_{j-c})$ such that

$\nu^{1-c}\t_c=-\sum_{i+i'=-c}\nu^iR\nu^{i'}$.
\nl
For $k\in[0,-c]$ we define $\t_{c,k}\in\Hom(P^\nu_c,\bV_{c+2k+j})$ by $\t_{c,0}=\t_c$ 
and $\t_{c,k}=\nu\t_{c,k-1}+R\nu^{k-1}$ for $k\in[1,-c]$. Then 
$\nu\t_{c,-c}+R\nu^{-c}=0$. Let $T:\bV@>>>\bV$ be the unique linear map such that 
$T(\nu^kx)=\t_{c,k}(x)$ for $x\in P^\nu_c,c\le0,k\in[0,-c]$. This $T$ has the required
property. 

\subhead 2.2\endsubhead
Let $\cc$ be the category whose objects are $\kk$-vector spaces of finite dimension;
morphisms are linear maps. Let $V\in\cc$. A collection of subspaces 
$V_*=(V_{\ge a})_{a\in\ZZ}$ of $V$ is said to be a {\it filtration} of $V$ if
$V_{\ge a+1}\sub V_{\ge a}$ for all $a$, $V_{\ge a}=0$ for some $a$, $V_{\ge a}=V$ for 
some $a$. We say that $V$ is {\it filtered} if a filtration $V_*$ of $V$ is given. 
Assume that this is the case. We set $\gr V_*=\op_{a\in\ZZ}\gr_aV_*\in\bcc$ where
$\gr_aV_*=V_{\ge a}/V_{\ge a+1}$. For any $j\in\ZZ$ let 
$E_{\ge j}V_*=\{T\in\End(V);T(V_{\ge a})\sub V_{\ge a+j}\qua\frl a\}$. Any such $T$ 
induces a linear map $\bT\in\End_j(\gr V_*)$.

\subhead 2.3\endsubhead
Let $V\in\cc$. Let $\Nil(V)=\{T\in\End(V);T\text{ nilpotent }\}$. Let $N\in\Nil(V)$. 
When $p=1$, the Dynkin-Kostant theory associates to $1+N$ a canonical filtration 
$V^N_*$ of $V$; in terms of a basis of $V$ of the form 

(a) $\{N^kv_r;r\in[1,t],k\in[0,e_r-1]\}$ with $v_r\in V,e_r\ge1,N^{e_r}v_r=0$ for 
$r\in[1,t]$,
\nl
$V^N_{\ge a}$ is the subspace spanned by 
$\{N^kv_r;r\in[1,t],k\in[0,e_r-1],2k+1\ge e_r+a\}$. This subspace makes sense for any 
$p$ and we denote it in general by $V^N_{\ge a}$; it is independent of the choice of 
basis: we have

$V^N_{\ge a}=\sum_{j\ge\max(0,a)}N^j(\ker N^{2j-a+1})$. 
\nl
The subspaces $V^N_{\ge a}$ form a filtration $V^N_*$ of $V$; thus, $V$ becomes a 
filtered vector space. From the definitions we see that

(b) {\it$N\in E_{\ge2}V^N_*$ and $\bN\in\End_2(\gr V_*^N)$ belongs to 
$\End_2^0(\gr V_*^N)$.}
\nl
Note that for any $j\ge1$, 

(c) {\it$\dim P^{\bN}_{1-j}$ is the number of Jordan blocks of size $j$ of $N:V@>>>V$.}
\nl
From 2.1(a) we deduce that for any $n\ge0$: 

(d) $\dim P^{\bN}_{-n}=\dim\gr_{-n}V^N_*-\dim\gr_{-n-2}V^N_*$.

\subhead 2.4\endsubhead
According to \cite{\DEII, 1.6.1},

(a) {\it if $V_*$ is a filtration of $V$ and $N\in E_{\ge2}V_*$ induces an element 
$\nu\in\End_2^0(\gr V_*)$ then $V_*=V_*^N$.}
\nl
We show that $V_{\ge a}=V^N_{\ge a}$ for all $a$. Let $e$ be the smallest integer 
$\ge0$ such that $N^e=0$. We argue by induction on $e$. If $a\ge e$ then 
$\nu^a:\gr_{-a}V_*@>>>\gr_aV_*$ is both $0$ and an isomorphism hence 
$V_{\ge-a}=V_{\ge1-a}$ and $V_{\ge a}=V_{\ge a+1}$. Thus, 
$V_{\ge e}=V_{\ge e+1}=\do=0$ and $V_{\ge1-e}=V_{\ge-e}=\do=V$. Similarly,
$V^N_{\ge e}=V^N_{\ge e+1}=\do=0$ and $V^N_{\ge1-e}=V^N_{\ge-e}=\do=V$. Hence 
$V_{\ge a}=V^N_{\ge a}$ if $a\ge e$ or if $a\le1-e$. This already suffices in the case 
where $e\le1$. Thus we may assume that $e\ge2$. Now
$\nu^{e-1}:\gr_{1-e}V_*@>>>\gr_{e-1}V_*$ is an isomorphism that is, 
$N^{e-1}:V/V_{\ge2-e}@>>>V_{\ge e-1}$ is an isomorphism. We see that 
$V_{\ge e-1}=N^{e-1}V$ and $V_{\ge2-e}=\ker(N^{e-1})$. Hence if $2-e\le a\le e-1$ we 
have $N^{e-1}V\sub V_{\ge a}\sub\ker(N^{e-1})$; let $V'_{\ge a}$ be the image of 
$V_{\ge a}$ under the obvious map $\r:\ker(N^{e-1})@>>>V':=\ker(N^{e-1})/N^{e-1}V$. For
$a\le1-e$ we set $V'_{\ge a}=V'$ and for $a\ge e$ we set $V'_{\ge a}=0$. Now 
$(V'_{\ge a})_{a\in\ZZ}$ is a filtration of $V'$ satisfying a property like (a) (with 
$N$ replaced by the map $N':V'@>>>V'$ induced by $N$). Since $N'{}^{e-1}=0$, the 
induction hypothesis applies to $N'$; it shows that $V'_{\ge a}=V'{}^{N'}_{\ge a}$ for 
all $a$. Since for $2-e\le a\le e-1$, $V_{\ge a}=\r\i(V'_{\ge a})$, it follows that 
$V_{\ge a}=\r\i(V'{}^{N'}_{\ge a})$; similarly, $V^N_{\ge a}=\r\i(V'{}^{N'}_{\ge a})$
hence $V_{\ge a}=V^N_{\ge a}$. This completes the proof.

With notation in the proof above we have:

$V^N_{\ge a}=0$ for $a\ge e$,

$V^N_{\ge a}=V$ for $a\le1-e$,

$V^N_{\ge a}=\r\i(V'{}^{N'}_{\ge a}),V'{}^{N'}_{\ge a}=\r(V^N_{\ge a})$ for $e\ge2$ and
$2-e\le a\le e-1$,

$V^N_{\ge e-1}=N^{e-1}V$ if $e\ge1$,

$V^N_{\ge2-e}=\ker(N^{e-1})$ if $e\ge1$.

We have $\gr_aV^N_*=0$ for $a\ge e$ and for $a\le-e$. 

Note also that the proof above provides an alternative (inductive) definition of
$V^N_{\ge a}$ which does not use a choice of basis.

\subhead 2.5\endsubhead
Let $V,N$ be as in 2.3. Let $V_*=V^N_*$. Let $\nu=\bN\in\End_2(\gr V_*)$. We can find a
grading $V=\op_{a\in\ZZ}V_a$ of $V$ such that 

(a) $NV_a\sub V_{a+2}$ and $V_{\ge a}=V_a\op V_{a+1}\op\do$ for all $a$. 
\nl
For example, in terms of a basis of $V$ as in 2.3(a), we can take $V_a$ to be the 
subspace spanned by $\{N^kv_r;r\in[1,t],k\in[0,e_r-1],2k+1=e_r+a\}$. Taking direct sum 
of the obvious isomorphisms $V_a@>\si>>\gr_aV_*$ we obtain an isomorphism of graded 
vector spaces $V@>\si>>\gr V_*$ under which $N$ corresponds to $\nu$. It follows that 

(b) $N\in\End_2^0(V)$ (defined in terms of the grading $\op_aV_a$).
\nl
We note the following result.

(c) {\it Let $n\ge0$ and let $x\in P^\nu_{-n}$. There exists a representative $\dx$ of 
$x$ in $V_{\ge-n}$ such that $N^{n+1}\dx=0$.}
\nl
Let $V_a$ be as above. There is a unique representative $\dx$ of $x$ in $V_{\ge-n}$ 
such that $\dx\in V_{-n}$. We have $N^{n+1}\dx\in V_n$ and the image of $N^{n+1}\dx$ 
under the canonical isomorphism $V_n@>\si>>\gr_nV_*^N$ is $0$; hence $N^{n+1}\dx=0$.

Let $E_{\ge1}^NV_*=\{S\in E_{\ge1}V_*;SN=NS\}$,
$\End_1^\nu(\gr V_*)=\{\s\in\End_1(\gr V_*),\s\nu=\nu\s\}$. We show:

(d) {\it The obvious map $E_{\ge1}^NV_*@>>>\End_1^\nu(\gr V_*),S\m\bS$ is surjective.}
\nl
Let $\s\in\End_1^\nu(\gr V_*)$. Let $V_a$ be as above. In terms of these $V_a$ we 
define $V@>\si>>\gr V_*$ as above. Under this isomorphism, $\s$ corresponds to a linear
map $S:V@>>>V$. Clearly, $S\in E_{\ge1}^NV_*$ and $\bS=\s$.

\subhead 2.6\endsubhead
Let $V,N$ be as in 2.3. Let $V_*=V_*^N$. Now $1+E_{\ge1}V_*$ is a subgroup of $GL(V)$
acting on $N+E_{\ge3}V_*$ by conjugation. We show that 

(a) {\it the conjugation action of $1+E_{\ge1}V_*$ on $N+E_{\ge3}V_*$ is transitive.}
\nl
We must show: if $S\in E_{\ge3}V_*$ then there exists $T\in E_{\ge1}V_*$ such that 
$(1+T)N=(N+S)(1+T)$ that is, $TN-NT=S+ST$. We fix subspaces $V_a$ as in 2.5. We have
$S=\sum_{j\ge3}S_j$ where $S_j\in\End(V)$ satisfy $S_jV_a\sub V_{a+j}$ for all $a$. We
seek a linear map $T=\sum_{j\ge1}T_j$ where $T_j\in\End(V)$ satisfy 
$T_jV_a\sub V_{a+j}$ for all $a$ and
$\sum_{j\ge1}(T_jN-NT_j)=\sum_{j\ge3}S_j+\sum_{j'\ge3,j''\ge1}S_{j'}T_{j''}$ that is,

$(*)$ $T_jN-NT_j=S_{j+2}+\sum_{j'\in[1,j-1]}S_{j+2-j'}T_{j'}$ for $j=1,2,\do$.   
\nl
We show that this system of equations in $T_j$ has a solution. We take $T_1=0$. Assume 
that $T_j$ has been found for $j<j_0$ for some $j_0\ge2$ so that $(*)$ holds for 
$j<j_0$. We set $R=S_{j_0+2}+\sum_{j'\in[1,j_0-1]}S_{j+2-j'}T_{j'}$. Then 
$R(V_a)\sub V_{a+j_0+2}$ for any $a$. The equation $T_{j_0}N-NT_{j_0}=R$ can be solved
by 2.1(b) (see 2.5(b)). This shows by induction that the system $(*)$ has a solution. 
(a) is proved.

We now show:

(b) {\it if $\tN\in N+E_{\ge3}V_*$ then $V^{\tN}_*=V_*$.}
\nl
Indeed by (a) we can find $u\in1+E_{\ge1}V_*$ such that $\tN=uNu\i$. Since $V_*^N$ is
canonically attached to $N$, we have $V^{uNu\i}_{\ge a}=u(V^N_{\ge a})=V^N_{\ge a}$ and
(b) follows. For example,

(c) {\it if $\tN=c_1N+c_2N^2+\do+c_kN^k$ where $c_i\in\kk,c_1\ne0$ then 
$V^{\tN}_*=V^N_*$.}
\nl
We may assume that $c_1=1$. Since $c_2N^2+\do+c_kN^k\in E_{\ge4}V_*\sub E_{\ge3}V_*$,
(b) is applicable and (c) follows.

\subhead 2.7\endsubhead
Let $V,N$ be as in 2.3. Let $V_*=V_*^N$. Let $\nu=\bN\in\End_2(\gr V_*)$. Let $r\ge2$ 
be such that $N^r=0$ on $V$. Let $W$ be an $N$-stable subspace of $V$ such that there 
exists an $N$-stable complement of $W$ in $V$, $N:W@>>>W$ has no Jordan block of size 
$\ne r,r-1$ and $N^{r-2}=0$ on $V/W$. Then $W_*=W_*^N$ is defined. Define a linear map 
$\mu:\gr V_*@>>>\gr W_*$ as follows. Let $x\in\gr_aV_*$. We have uniquely 
$x=\sum_{c\in a+2\ZZ;c\le a\le-c}\nu^{(a-c)/2}x_c$ where $x_c\in P^\nu_c$; we set
$$\mu(x)=\sum_{c\in a+2\ZZ;c\le a\le-c,c=1-r\text{ or }2-r}\nu^{(a-c)/2}x_c.$$
Let $\cx$ be the set of $N$-stable complements of $W$ in $V$. Then $\cx\ne\em$. For 
$Z\in\cx$ define $\Pi_Z:V@>>>W$ by $\Pi_Z(w+z)=w$ where $w\in W,z\in Z$. Let 
$\baP_Z:\gr V_*@>>>\gr W_*$ be the map induced by $\Pi_Z$. We show that 

(a) $\Pi_Z(V_{\ge a})\sub W_{\ge a}$ for all $a$ and $\baP_Z=\mu$. 
\nl
We have $V_{\ge a}=W_{\ge a}\op Z_{\ge a}$. If $x\in V_{\ge a},x=w+z,w\in W_{\ge a},
z\in Z_{\ge a}$, then $\Pi_Z(x)=w$. Thus $\Pi_Z(V_{\ge a})\sub W_{\ge a}$. We can find 
direct sum decompositions $W=\op_mW_m,Z=\op_mZ_m$ such that $NW_m\sub W_{m+2}$,
$NZ_m\sub Z_{m+2}$ and $N^m:W_{-m}@>\si>>W_m$, $N^m:Z_{-m}@>\si>>Z_m$ for $m\ge0$ (see 
2.5). Let $V_a=W_a\op Z_a$. Define $V_a^{\prim},W_a^{\prim},Z_a^{\prim}$ as in 2.1 in 
terms of $N$. We have $V_a^{\prim}=W_a^{\prim}\op Z_a^{\prim}$. We must show that 
$\baP_Z(x)=\mu(x)$ for $x\in\gr_aV_*$. It suffices to show: if $w\in W_a,z\in Z_a$ and
$w+z=\sum_{c\in a+2\ZZ;c\le a\le-c}\nu^{(a-c)/2}x_c$ where $x_c\in V_c^{\prim}$ then 
$w=\sum_{c\in a+2\ZZ;c\le a\le-c,1-r\le c\le2-r}\nu^{(a-c)/2}x_c$. We have 
$x_c=w_c+z_c$ where $w_c\in W_c^{\prim},z_c\in Z_c^{\prim}$ and 
$w=\sum_{c\in a+2\ZZ;c\le a\le-c}\nu^{(a-c)/2}w_c$. Now if $W_c^{\prim}\ne0$ then 
$1-r\le c\le2-r$. Hence 
$w=\sum_{c\in a+2\ZZ;c\le a\le-c,1-r\le c\le2-r}\nu^{(a-c)/2}w_c$.

Also, $Z_{1-r}^{\prim}=Z_{2-r}^{\prim}=0$ since $N:Z@>>>Z$ has no Jordan blocks of size
$\ge r-1$. Thus if $c\in a+2\ZZ;c\le a\le-c,1-r\le c\le2-r$ then $z_c=0$ and 
$x_c=w_c$. Thus $w=\sum_{c\in a+2\ZZ;c\le a\le-c,1-r\le c\le2-r}\nu^{(a-c)/2}x_c$, as 
required.

Let $Z,Z'\in\cx$. By the previous argument, $\Pi_Z,\Pi_{Z'}:V@>>>W$ both map 
$V_{\ge a}$ into $W_{\ge a}$ and induce the same map $\gr V_*@>>>\gr W_*$. It follows 
that $\Pi_Z-\Pi_{Z'}:V@>>>W$ maps $V_{\ge a}$ into $W_{\ge a+1}$. In other words, 

(b) {\it if $x\in V_{\ge a}$ and $x=w+z=w'+z'$ where $w,w'\in W,z\in Z,z'\in Z'$, then 
$w-w'\in W_{\ge a+1}$.}
\nl
Define $\Ph\in GL(V)$ by $\Ph(x)=x$ for $x\in W$, $\Ph(x)=x'$ for $x\in Z$ where 
$x'\in Z'$ is given by $x-x'\in W$. We show:

(c) $(1-\Ph)V_{\ge a}\sub V_{\ge a+1}$ {\it for any $a$.}
\nl
Let $x\in V_{\ge a}$. We have $x=w+z=w'+z'$ where $w,w'\in W,z\in Z,z'\in Z'$. We have 
$\Ph(x)=w+z'$ hence $(1-\Ph)(x)=(w+z)-(w+z')=z-z'=w'-w$ and this belongs to 
$W_{\ge a+1}$ by (b).

We show: 

(d) $\Ph N=N\Ph$.
\nl
Indeed, for $x=x_1+x_2,x_1\in W,x_2\in Z$ we have $Nx=Nx_1+Nx_2$ with $Nx_1\in W$,
$Nx_2\in Z$ and $x_2-x'_2\in W$ with $x'_2\in Z'$. We have $Nx_2-Nx'_2\in W$ with 
$Nx_2\in Z$, $Nx'_2\in Z'$. Hence $\Ph(Nx)=Nx_1+Nx'_2=N(x_1+x'_2)=N\Ph(x)$, as 
required.

\subhead 2.8\endsubhead
Let $V,N$ be as in 2.3. Let $r\ge1$ be such that $N^r=0$. A subspace $W$ of $V$ is said
to be {\it $r$-special} if $NW\sub W$, $N:W@>>>W$ has no Jordan blocks of size $\ne r$ 
and $N^{r-1}=0$ on $N/W$. We show:

(a) {\it If $W,W'$ are $r$-special subspaces then there exists a subspace $X$ of $V$ 
such that $NX\sub X,W\op X=V,W'\op X=V$.}
\nl
We argue by induction on $r$. If $r=1$ the result is obvious; we have $W=W'=V$. Assume 
that $r\ge2$. Let $V'=\ker N^{r-1},V''=\ker N^{r-2}$. Let $E\sub W,E'\sub W'$ be such 
that $W=E\op NE\op\do N^{r-1}E$, $W'=E'\op NE'\op\do N^{r-1}E'$. Clearly,
$E\cap V'=0$, $E'\cap V'=0$, $NE\sub V'$, $NE\cap V''=0$, $NE'\sub V'$, 
$NE'\cap V''=0$. Let $E''$ be a subspace of $V'$ such that $E''$ is a complement of 
$NE\op V''$ in $V'$ and a complement of $NE'\op V''$ in $V'$. (Such $E''$ exists since 
$\dim(NE\op V'')=\dim(NE'\op V')=\dim E+\dim V''=\dim E'+\dim V''$.) Then 

$W_1=(E''\op NE)+N(E''\op NE)+\do+N^{r-2}(E''\op NE)$,

$W'_1=(E''\op NE')+N(E''\op NE')+\do+N^{r-2}(E''\op NE')$
\nl
are $(r-1)$-special subspaces of $V'$. By the induction hypothesis we can find an 
$N$-stable subspace $X_1$ of $V'$ such that $V_1\op X_1=V',V'_1\op X_1=V'$. Then 
$X=(E''+N(E'')+\do+N^{r-2}(E''))+X_1$ has the required properties.

(b) {\it If $W,W'$ are $r$-special subspaces then there exists $g\in1+E_{\ge1}V_*$ such
that $g(W)=W'$, $gN=Ng$.}
\nl
Let $X$ be as in (a). Define $g\in GL(V)$ by $g(x)=x$ for $x\in X$ and $g(w)=w'$ for 
$w\in W$ where $w'\in W'$ is given by $w-w'\in X$. Then $g(W)=W'$, $(g-1)X=0$ and
$(g-1)W\sub X$. Clearly, $gN=Ng$. We have $V_{\ge a}=W_{\ge a}\op X_{\ge a}$. It 
suffices to show that $(g-1)(W_{\ge a})\sub X_{\ge a+1}$. Now $X=X_{\ge2-r}$. We have
$W=W_{\ge1-r},W_{\ge2-r}=W_{\ge3-r}=NW,W_{\ge4-r}=W_{\ge5-r}=N^2W,\do$. Now if 
$a\le1-r$ then $(g-1)W_{\ge a}=(g-1)W\sub X=X_{\ge a+1}$. If $a=2-r$ or $a=3-r$ then 

$(g-1)W_{\ge a}=(g-1)NW=N(g-1)W\sub NX=NX_{\ge2-r}\sub X_{\ge4-r}\sub X_{\ge a+1}$.
\nl
If $a=4-r$ or $a=5-r$ then 

$(g-1)W_{\ge a}=(g-1)N^2W=N^2(g-1)W\sub N^2X=N^2X_{\ge2-r}\sub X_{\ge6-r}\sub 
X_{\ge a+1}$.
\nl
Continuing in this way, the result follows. 

\subhead 2.9\endsubhead
Let $V\in\bcc$. Let $G=GL(V)$. For any filtration $V_*$ of $V$ let 
$$\x(V_*)=\{N\in\Nil(V);V^N_*=V_*\}=\{N\in E_{\ge2}V_*;\bN\in\End_2^0(\gr V_*)\}$$ 
(see 2.3(b), 2.4(a)). The following three conditions are equivalent:

(i) $\x(V_*)\ne\em$;

(ii) $\End_2^0(\gr V_*)\ne\em$;

(iii) $\dim\gr_nV_*=\dim\gr_{-n}V_*\ge\dim\gr_{-n-2}V_*$ for all $n\ge0$.
\nl
We have (i)$\imp$(ii) by the definition of $\x(V_*)$; we have (ii)$\imp$(iii) by
2.3(d). The fact that (iii)$\imp$(ii) is easily checked. If (ii) holds we pick for any
$a$ a subspace $V_a$ of $V_{\ge a}$ complementary to $V_{\ge a+1}$ and an element in 
$\End_2^0(V)$ (defined in terms of the grading $\op_aV_a$). This element is in 
$\x(V_*)$ and (i) holds. 

Let $\fF_V$ be the set of all filtrations $V_*$ of $V$ that satisfy (i)-(iii). From the
definitions we have a bijection 

(a) $\fF_V@>\si>>D_G,V_*\m\l$
\nl
($D_G$ as in 1.1) where $\l=(G^\l_0\sps G^\l_1\sps G^\l_2\sps\do)$ is defined in terms 
of $V_*$ by $G^\l_0=E_{\ge0}V_*\cap G$ and $G^\l_n=1+E_{\ge n}V_*$ for $n\ge1$.

The sets $\x(V_*)$ (with $V_*\in\fF_V)$ form a partition of $\Nil(V)$. (If 
$N\in\Nil(V)$ we have $N\in\x(V_*)$ where $V_*=V^N_*$).

Let $V_*\in\fF_V$. Let $\Pi=E_{\ge0}V_*\cap G$. We show that $\x(V_*)$ is a single 
$\Pi$-conjugacy class. Let $N,N'\in\x(V_*)$. Since $V^N_*=V^{N'}_*$ we see from 2.3(d) 
that $\dim P^{\bN}_{1-j}=\dim P^{\bN'}_{1-j}$ for any $j\ge0$. Using 2.3(c) we see that
for any $j\ge0$, $N,N'$ have the same number of Jordan blocks of size $j$. Hence there 
exists $g\in G$ such that $N'=gNg\i$. For any $a$, 
$gV_{\ge a}^N=V_{\ge a}^{N'}=V_{\ge a}^N$ hence $gV_{\ge a}=V_{\ge a}$. We see that 
$g\in E_{\ge0}$ hence $g\in\Pi$, as required. Taking in the previous argument $N'=N$, 
we see that, if $N\in\x(V_*)$ and $g\in G$ satisfies $gNg\i=N$ then $g\in\Pi$. Now any 
element in $E_{\ge2}V_*-\x(V_*)$ is in the closure of $\x(V_*)$ (since $E_{\ge2}V_*$ is
irreducible and $\x(V_*)$ is open in it (and non-empty) hence it is in the closure of 
the $G$-conjugacy class containing $\x(V_*)$. We show that it is not contained in that 
$G$-conjugacy class. (Assume that it is. Then we can find $N\in\x(V_*)$ and 
$N'\in E_{\ge2}V_*-\x(V_*)$ that are $G$-conjugate. Then the $\Pi$-orbit $\Pi(N)$ of 
$N$ in $E_{\ge2}V_*$ is $\x(V_*)$ hence is dense in $E_{\ge2}V_*$ while the $\Pi$-orbit
$\Pi(N')$ of $N'$ is contained in the proper closed subset $E_{\ge2}V_*-\x(V_*)$ of 
$E_{\ge2}V_*$; hence $\dim\Pi(N)=\dim(E_{\ge2}V_*)>\dim\Pi(N')$. It follows that $a<a'$
where $a$ (resp. $a'$) is the dimension of the centralizer of $N$ (resp. $N'$) in 
$\Pi$. Let $\ta$ (resp. $\ta'$) be the dimension of the centralizer of $N$ (resp. $N'$)
in $G$. By an earlier argument we have $a=\ta$. Obviously, $a'\le\ta'$. Since $N,N'$ 
are $G$-conjugate, we have $\ta=\ta'$. Thus, $\ta=a<a'\le\ta'=\ta$, contradiction.) We 
see that $1+\x(V^*)=X^\l$ where $V_*\m\l$ as in (a) and $X^\l$ is as in 1.1. Thus 
$\fP_1$ holds for $G$. From this $\fP_2,\fP_3$ follow; $H^\bla$ in $\fP_2$ is a single 
conjugacy class in this case. Also, $\fP_8$ is trivial since $G^\l_0$ acts transitively
on $X^\l$. Now $\fP_5$ is easily verified. $\fP_6$ (hence $\fP_4$) follows from 2.6(a);
$\fP_7$ is trivial in this case.

\head 3. Symplectic groups\endhead
\subhead 3.1\endsubhead
In this section, any text marked as $\sp\do\sp$ applies only in the case $p=2$. 

For $V,V'\in\cc$ let $\Bil(V,V')$ be the space of all bilinear forms $V\T V'@>>>\kk$. 
For $b\in\Bil(V,V')$ define $b^*\in\Bil(V',V)$ by $b^*(x,y)=b(y,x)$. We write $\Bil(V)$
instead of $\Bil(V,V)$. Let $\Symp(V)$ be the set of non-degenerate symplectic forms on
$V$.

Let $\bV\in\bcc$. We say that $\lar_0\in\Symp(\bV)$ is {\it admissible} if 
$\la x,y\ra_0=0$ for $x\in\bV_a,y\in\bV_{a'},a+a'\ne0$. Assume that 
$\lar_0\in\Symp(\bV)$ is admissible and that $\nu\in\End_2^0(\bV)$ is {\it 
skew-adjoint} that is, $\la\nu(x),y\ra_0+\la x,\nu(y)\ra_0=0$ for $x,y\in\bV$. For 
$n\ge0$ we define a bilinear form $b_n:P^\nu_{-n}\T P^\nu_{-n}@>>>\kk$ by 
$b_n(x,y)=\la x,\nu^ny\ra_0$. We show: 

(a) $b_n(x,y)=(-1)^{n+1}b_n(y,x)$ for $x,y\in P^\nu_{-n}$.
\nl
Indeed,
$$b_n(x,y)=\la x,\nu^ny\ra_0=(-1)^n\la\nu^nx,y\ra_0=(-1)^{n+1}\la y,\nu^nx\ra_0
=(-1)^{n+1}b_n(y,x),$$
as required. We show:

(b) {\it$b_n$ is non-degenerate.}
\nl
Let $y\in P^\nu_{-n}$ be such that $\la x,\nu^n y\ra_0=0$ for all $x\in P^\nu_{-n}$. If
$x'\in P^\nu_{-n-2k},k>0$, we have 
$\la\nu^kx',\nu^ny\ra_0=\pm\la x,\nu^{n+k}y\ra_0=\pm\la x,0\ra_0=0$. Since 
$\bV_{-n}=\sum_{k\ge0}\nu^kP^\nu_{-n-2k}$, we see that $\la x,\nu^ny\ra_0=0$ for all 
$x\in\bV_{-n}$. Since $\la\bV_m,\nu^ny\ra_0=0$ for $m\ne-n$ we see that 
$\la\bV,\nu^ny\ra_0=0$. By the non-degeneracy of $\lar_0$, it follows that $\nu^ny=0$. 
Since $\nu^n:\bV_{-n}@>\si>>\bV_n$, it follows that $y=0$, as required. We show:

(c) {\it if $n\ge0$ is even then $b_n$ is a symplectic form. Hence $\dim P^\nu_{-n}$ is
even.}
\nl
Indeed, for $x\in P^\nu_{-n}$ we have 
$\la x,\nu^nx\ra_0=\pm\la\nu^{n/2}x,\nu^{n/2}x\ra_0=0$.

\subhead 3.2\endsubhead
Let $V\in\cc$ and let $\lar\in\Symp(V)$. Let 

$Sp(\lar)=\{T\in GL(V);T\text{ preserves }\lar\}$.
\nl
For any subspace $W$ of $V$ we set $W^\pe=\{x\in V;\la x,W\ra=0\}$. A filtration $V_*$ 
of $V$ is said to be {\it self-dual} if $(V_{\ge a})^\pe=V_{\ge1-a}$ for any $a$. It 
follows that 

(a) $\la V_{\ge a},V_{\ge a'}\ra=0$ if $a+a'\ge1$.
\nl
It also follows that there is a unique admissible $\lar_0\in\Symp(\gr V_*)$ such that 
for $x\in\gr_aV_*,y\in\gr_{-a}V_*$ we have $\la x,y\ra_0=\la\dx,\dy\ra$ where 
$\dx\in V_{\ge a},\dy\in V_{\ge-a}$ represent $x,y$. Moreover, 

(b) {\it there exists a direct sum decomposition $\op_{a\in\ZZ}V_a$ of $V$ such that
$V_{\ge a}=V_a\op V_{a+1}\op\do$ for all $a$ and $\la V_a,V_{a'}\ra=0$ for all $a,a'$ 
such that $a+a'\ne0$.}

Let $\cm_{\lar}$ be the set of $N\in\Nil(V)$ such that 
$\la Nx,y\ra+\la x,Ny\ra+\la Nx,Ny\ra=0$ for all $x,y$ or equivalently 
$1+N\in Sp(\lar)$. Define an involution $N\m N^\da$ of $\cm_{\lar}$ by 
$\la x,Ny\ra=\la N^\da x,y\ra$ for all $x,y\in V$ or equivalently by 
$N^\da=(1+N)\i-1=-N+N^2-N^3+\do$.

Let $N\in\cm_{\lar}$. We set $V_*=V^N_*$. By 2.6(c) we have $V^{N^\da}_*=V_*$. We 
show:

(c) {\it the filtration $V_*$ is self-dual.}
\nl
We argue by induction on $e$ as in 2.4. If $a\ge e$ then $V_{\ge a}=0,V_{\ge1-a}=V$ and
(c) holds. If $a\le1-e$ then $V_{\ge a}=V,V_{\ge1-a}=0$ and (c) holds. If $e\le1$ this 
already suffices. Hence we may assume that $e\ge2$ and $2-e\le a\le e-1$ hence 
$2-e\le1-a\le e-1$. Let $V'=\ker(N^{e-1})/\Im(N^{e-1})$. Let $\r:\ker(N^{e-1})@>>>V'$ 
be the canonical map. We have $N^{e-1}V=\ker((N^\da)^{e-1})^\pe=\ker(N^{e-1})^\pe$ 
since $(N^\da)^{e-1}=(-N)^{e-1}$. Hence $\lar$ induces $\lar'\in\Symp(V')$. Also $N$ 
induces a linear map $N':V'@>>>V'$ such that $N'\in\cm_{\lar'}$. By the induction 
hypothesis, $V'{}^{N'}_{\ge1-a}$ is the perpendicular in $V'$ of $V'{}^{N'}_{\ge a}$. 
Hence $V_{\ge a}=\r\i(V'{}^{N'}_{\ge a})$ is the perpendicular in 
$V$ of $V_{\ge1-a}=\r\i(V'{}^{N'}_{\ge1-a})$. This completes the proof.

Let $\nu\in\End_2^0(\gr V_*)$ be the endomorphism induced by $N$. We show that 

(d) {\it$\nu$ is skew-adjoint (with respect to $\lar_0$ on $\gr V^N_*$).}
\nl
It suffices to show that, if $a+a'+2=0$ and $x\in V_{\ge a'}$, $y\in V_{\ge a}$ then 
$\la Nx,y\ra+\la x,Ny\ra=0$. It suffices to show that $\la Nx,Ny\ra=0$. From (a),(b) we
see that $\la V_{\ge-1-a},Ny\ra=0$ hence it suffices to show that $Nx\in V_{\ge-1-a}$. 
We have $Nx\in V_{\ge a'+2}\sub V_{\ge-1-a}$ since $a'+2>-1-a$. This proves (c).

\subhead 3.3\endsubhead
$\sp$ In this subsection we assume that $p=2$. Let $V,\lar,N,\nu,\lar_0$ be as in 3.2. 
Let $V_*=V^N_*$. Then $b_n\in\Bil(P_{-n}^\nu)$ is defined for $n\ge0$, see 3.1. Let 
$\cl$ be the set of all even integers $n\ge2$ such that $b_{n-1},b_{n+1}$ are 
symplectic forms. Let $\cl'$ be the set of all even integers $n\ge2$ such that 
$b_{n-1},b_{n+1},b_{n+3},\do$ are symplectic forms or equivalently, if
$\la z,\nu^{n-1}(z)\ra_0=0$ for all $z\in\gr_{1-n}V_*$. (Assume first that 
$b_{n-1},b_{n+1},b_{n+3},\do$ are symplectic forms. By 2.1(a), any $z\in\gr_{1-n}V_*$ 
is of the form $\sum_{k\ge0}\nu^kz_k$ where $z_k\in P^\nu_{1-n-2k}$. For $k\ge0$ we 
have $\la\nu^kz_k,\nu^{n-1}(\nu^kz_k)\ra_0=0$ since $b_{n+2k-1}$ is symplectic. Since 
$z'\m\la z',\nu^{n-1}(z')\ra_0$ is additive in $z'$ it follows that 
$\la z,\nu^{n-1}(z)\ra_0=0$. Conversely, assume that $\la z,\nu^{n-1}(z)\ra_0=0$ for 
any $z\in\gr_{1-n}V_*$. In particular, for $k\ge0$ and $z_k\in P^\nu_{1-n-2k}$ we have 
$\la \nu^kz_k,\nu^{n-1}(\nu^kz_k)\ra_0=0$ that is, $\la z_k,\nu^{n-1+2k}z_k)\ra_0=0$. 
We see that $b_{n+2k-1}$ is symplectic.)

Clearly, $\cl'\sub\cl$. 

For $n\in\cl$, we define $q_n:P^\nu_{-n}@>>>\kk$ by $q_n(x)=\la\dx,N^{n-1}\dx\ra$ where
$\dx\in V_{\ge-n}$ is a representative for $x\in P^\nu_{-n}$ such that $N^{n+1}\dx=0$ 
(see 2.5(c)). We show that $q_n(x)$ is well defined. It suffices to show that if 
$y\in V_{\ge1-n},N^{n+1}y=0$ then $\la\dx+y,N^{n-1}(\dx+y)\ra=\la\dx,N^{n-1}\dx\ra$ 
that is, $\la y,N^{n-1}(y)\ra+\la\dx,N^{n-1}(y)\ra+\la y,N^{n-1}(\dx)\ra=0$. Since
$N^{n+1}(\dx)=0$, we have 
$$\la\dx,N^{n-1}(y)\ra+\la y,N^{n-1}(\dx)\ra=\la y,(N^{n-1}+(N^\da)^{n-1})(\dx)\ra
=\la y,N^n(\dx)\ra.$$
This is zero, since $y\in V_{\ge1-n},N^n(\dx)\in V_{\ge n}$ and $1-n+n=1$. It remains 
to show that $\la y,N^{n-1}(y)\ra=0$. It suffices to show that 
$\la z,\nu^{n-1}(z)\ra_0=0$ for all $z\in\gr_{1-n}V_*$ such that $N^{n+1}z=0$. By 
2.1(a) any such $z$ is of the form $z_0+\nu z_1$ where 
$z_0\in P^\nu_{1-n},z_1\in P^\nu_{-1-n}$. Now $z'\m\la z',\nu^{n-1}(z')\ra_0$ is 
additive in $z'$ hence it suffices to show that $\la z_0,\nu^{n-1}(z_0)\ra_0=0$ and 
$\la\nu(z_1),\nu^{n-1}(\nu(z_1))\ra_0=0$ for $z_0,z_1$ as above. This follows from our 
assumption that $b_{n-1}$ and $b_{n+1}$ are symplectic.

We show:

(a) {\it For $x,y\in P^\nu_{-n}$ we have $q_n(x+y)=q_n(x)+q_n(y)+b_n(x,y)$.}
\nl
Let $\dx,\dy\in V_{\ge-n}$ be representatives for $x,y$ such that $N^{n+1}\dx=0$,
$N^{n+1}\dy=0$. We must show that

$\la\dx+\dy,N^{n-1}(\dx+\dy)\ra=\la\dx,N^{n-1}(\dx)\ra+\la\dy,N^{n-1}(\dy)\ra
+\la\dx,N^n(\dy)\ra$,
\nl
or that 

$\la\dx,N^{n-1}(\dy)\ra+\la\dy,N^{n-1}(\dx)\ra+\la\dx,N^n(\dy)\ra=0$,
\nl
or that $\la\dx,((N^\da)^{n-1}+N^{n-1}+N^n)\dy\ra=0$. Since $n$ is even, 
$(N^\da)^{n-1}+N^{n-1}+N^n$ is a linear combination of $N^{n+1},N^{n+2},\do$ and it
remains to use that $N^{n+1}(\dy)=0$.

For $n\in\cl'$, we define $Q_n:\gr_{-n}V_*@>>>\kk$ by $Q_n(x)=\la\dx,N^{n-1}\dx\ra$ 
where $\dx\in V_{\ge-n}$ is a representative for $x$. We show that $Q_n(x)$ is well 
defined. It suffices to show that if $y\in V_{\ge1-n}$ then 
$\la\dx+y,N^{n-1}(\dx+y)\ra=\la\dx,N^{n-1}\dx\ra$ that is,
$\la y,N^{n-1}(y)\ra+\la\dx,N^{n-1}(y)\ra+\la y,N^{n-1}(\dx)\ra=0$. We have 

$\la\dx,N^{n-1}(y)\ra+\la y,N^{n-1}(\dx)\ra=\la y,(N^{n-1}+(N^\da)^{n-1})(\dx)\ra$
\nl
and this is a linear combination of terms $\la y,N^{n'}(\dx)\ra$ with $n'\ge n$. Each 
of these terms is $0$ since $y\in V_{\ge1-n},N^{n'}(\dx)\in V_{\ge2n'-n}$ and
$1-n+2n'-n\ge1$. It remains to show that $\la y,N^{n-1}(y)\ra=0$. This follows from the
fact that $\la z,\nu^{n-1}(z)\ra_0=0$ for all $z\in\gr_{1-n}V_*$. 

For $n\in\cl'$ we show:

(b) {\it if $x,y\in\gr_{-n}V_*$ then $Q_n(x+y)=Q_n(x)+Q_n(y)+\la x,\nu^ny\ra$.}
\nl
Let $\dx,\dy\in V_{\ge-n}$ be representatives for $x,y$. We must show that

$\la\dx+\dy,N^{n-1}(\dx+\dy)\ra=\la\dx,N^{n-1}(\dx)\ra+\la\dy,N^{n-1}(\dy)\ra
+\la\dx,N^n(\dy)\ra$,
\nl
or that 

$\la\dx,N^{n-1}(\dy)\ra+\la\dy,N^{n-1}(\dx)\ra+\la\dx,N^n(\dy)\ra=0$, 
\nl
or that $\la\dx,((N^\da)^{n-1}+N^{n-1}+N^n)\dy\ra$ is $0$. Since $n$ is even, this is a
linear combination of terms $\la \dx,N^{n'}(\dy)\ra$ with $n'>n$. Each of these terms 
is $0$ since $N^{n'}(\dy)\in V_{\ge2n'-n},\dx\in V_{\ge-n}$ and $2n'-n-n\ge1$. 

Now let $n\in\cl'$ and let $x\in\gr_{-n}V_*$. We can write $x=\sum_{k\ge0}\nu^kx_k$ 
where $x_k\in P^\nu_{-n-2k}$. We show that

(c) $Q_n(x)=\sum_{k\ge0}q_{n+2k}(x_k)$.
\nl
Let $\dx_k$ be a representative of $x_k$ in $V_{\ge-n-2k}$ such that
$N^{n+2k+1}\dx_k=0$. Then $\sum_{k\ge0}N^k\dx_k$ is a representative of $x$ in 
$V_{\ge-n}$ and we must show:
$$\la\sum_{k\ge0}N^k\dx_k,N^{n-1}\sum_{k'\ge0}N^{k'}\dx_{k'}\ra=
\sum_{k\ge0}\la\dx_k,N^{n+2k-1}\dx_k\ra.$$
The left hand side is $\sum_{k,k'\ge0}\la N^k\dx_k,N^{n-1+k'}\dx_{k'}\ra$. If 
$k\ge k'+2$ we have 

$\la N^k\dx_k,N^{n-1+k'}\dx_{k'}\ra=\la\dx_k,(N^\da)^kN^{n-1+k'}\dx_{k'}\ra$
\nl
and this is zero since $N^{n+2k'+1}\dx_{k'}=0$. If $k'\ge k+2$ we have

$\la N^k\dx_k,N^{n-1+k'}\dx_{k'}\ra=\la (N^\da)^{n-1+k'}N^k\dx_k,\dx_{k'}\ra$
\nl
and this is zero since $N^{n+2k+1}\dx_k=0$. It suffices to show that 
$$\align&\sum_{k\ge0}(\la N^k\dx_k,N^{n-1+k}\dx_k\ra
+\la N^{k+1}\dx_{k+1},N^{n-1+k}\dx_k\ra+\la N^k\dx_k,N^{n+k}\dx_{k+1}\ra)\\&
=\sum_{k\ge0}\la \dx_k,N^{n+2k-1}\dx_k\ra.\endalign$$
We have
$$\align&\la N^{k+1}\dx_{k+1},N^{n-1+k}\dx_k\ra+\la N^k\dx_k,N^{n+k}\dx_{k+1}\ra\\&
=\la N^{k+1}\dx_{k+1},(N^{n-1+k}+(N^\da)^{n-1}N^k)\dx_k\ra\\&=\la N^{k+1}\dx_{k+1},
(c_1N^{n+k}+c_2N^{n+k+1}+\do)\dx_k\ra\\&=c_1\la \nu^{k+1}x_{k+1},
\nu^{n+k}x_k\ra_0=c_1\la x_{k+1},\nu^{n+2k+1}x_k\ra_0=0.\endalign$$
(Here $c_1,c_2,\do\in\kk$.) It suffices to show that 
$\la N^k\dx_k,N^{n-1+k}\dx_k\ra+\la\dx_k,N^{n+2k-1}\dx_k\ra$ is $0$. This equals  
$$\align&\la \dx_k,(N^{n+2k-1}+(N^\da)^kN^{n-1+k})\dx_k\ra=\la \dx_k,(N^{n+2k}
+c'_1N^{n+2k+1}+\do)\dx_k\ra\\&
=\la x_k,\nu^{n+2k}x_k\ra_0=\la \nu^{k+n/2}x_k,\nu^{k+n/2}x_k\ra_0=0.\endalign$$
(Here $c'_1,c'_2,\do\in\kk$.) This completes the proof of (c). 

We say that $(q_n)_{n\in\cl}$ are {\it the quadratic forms attached to $(N,\lar)$.} We 
say that $(Q_n)_{n\in\cl'}$ are {\it the Quadratic forms attached to $(N,\lar)$}. $\sp$

\subhead 3.4\endsubhead
Let $V\in\cc$ and let $V_*$ be a filtration of $V$. We fix 
$\lar_0\in\Symp(\gr V_*)$ which is admissible and $\nu\in\End_2^0(\gr V_*)$ which is
skew-adjoint with respect to $\lar_0$ (see 3.1). Then $P^\nu_{-n}$ are defined in terms
of $\gr V_*,\nu$ and $b_n\in\Bil(P^\nu_{-n})$ are defined as in 3.1 for any $n\ge0$. 
Let $\cv=1+E_{\ge1}V_*$, a subgroup of $GL(V)$. 

$\sp$ If $p=2$, let $\nn$ be the smallest even integer $\ge2$ such that  
$b_{\nn-1},b_{\nn+1},b_{\nn+3},\do$ are symplectic or, equivalently, such that
$\la z,\nu^{\nn-1}(z)\ra_0=0$ for all $z\in\gr_{1-\nn}V_*$. Let 
$Q:\gr_{-\nn}V_*@>>>\kk$ be a quadratic form such that 
$Q(x+y)=Q(x)+Q(y)+\la x,\nu^\nn y\ra$ for all $x,y\in\gr_{-\nn}V_*$. $\sp$.

Let $\cz$ be the set of all pairs $(N,\lar)$ where $N\in\Nil(V)$, $\lar\in\Symp(V)$ are
such that $V^N_*=V_*$, $\la Nx,y\ra+\la x,Ny\ra+\la Nx,Ny\ra=0$ for $x,y\in V$, $N$ 
induces $\nu$ on $\gr V_*$, $\lar$ induces $\lar_0$ on $\gr V_*$; $\sp$ in the case 
$p=2$ we require in addition that $Q_\nn$ defined in terms of $(N,\lar)$ as in 3.3 is 
equal to $Q$. $\sp$

The proofs of Propositions 3.5, 3.6, 3.7 below are intertwined (see 3.11).

\proclaim{Proposition 3.5}In the setup of 3.4 let $\lar\in\Symp(V)$ be such that $V_*$ 
is self-dual with respect to $\lar$ and $\lar$ induces $\lar_0$ on $\gr V_*$. Let 
$Y=Y_{\lar}=\{N;(N,\lar)\in\cz\}$. Let $U'=\cv\cap Sp(\lar)$, a subgroup of $Sp(\lar)$.
Then

(a) $Y\ne\em$;

(b) if $N\in Y$ and $z\in U'$ then $zNz\i\in Y$ (thus $U'$ acts an $Y$ by conjugation);

(c) the action (b) of $U'$ on $Y$ is transitive.
\endproclaim
The proof of (a) is given in 3.8. Now (b) follows immediately from 3.7(a).

We show that (c) is a consequence of 3.7(c). Assume that 3.7(c) holds. Let $N,N'\in Y$.
We have $(N,\lar)\in\cz,(N',\lar)\in\cz$ and by 3.7(c) there exists $g\in\cv$ such that
$N'=gNg\i$, $\la g\i x,g\i y\ra=\la x,y\ra$ for $x,y\in V$. Then $g\in U'$ and (c) is 
proved (assuming 3.7(c)).

\proclaim{Proposition 3.6}In the setup of 3.4 let $N\in\Nil(V)$ be such that 
$V^N_*=V_*$ and $N$ induces $\nu$ on $\gr V_*$. Let 
$X=X_N=\{\lar;(N,\lar)\in\cz\}$. Let $U=U_N=\{T\in\cv;TN=NT\}$, a subgroup of 
$GL(V)$. Then:

(a) $X\ne\em$;

(b) if $\lar\in X$ and $u\in U$ then the symplectic form $\lar'$ on $V$ given by 
$\la x,y\ra'=\la u\i x,u\i y\ra$ belongs to $X$ (thus $U$ acts naturally an $X$);

(c) the action (b) of $U$ on $X$ is transitive.
\endproclaim
We show that (a) is a consequence of 3.7(a). By 3.7(a) there exists 
$(N',\lar')\in\cz$. By 2.6(a) there exists $g\in\cv$ such that $N=gN'g\i$. Define 
$\lar\in\Symp(V)$ by $\lar=\la g\i x,g\i y\ra'$. From 3.7(a) we see that 
$(N,\lar)\in\cz$ hence $\lar\in X_N$. Thus $X_N\ne\em$, as required.

Now (b) follows immediately from 3.7(b). The proof of (c) is given in 3.9, 3.10.

\proclaim{Proposition 3.7}In the setup of 3.4,

(a) $\cz\ne\em$;

(b) if $(N,\lar)\in\cz$, $g\in\cv$ and $(N',\lar')$ is defined by $N'=gNg\i$,
$\la x,y\ra'=\la g\i x,g\i y\ra$ then $(N',\lar')\in\cz$ (thus $\cv$ acts naturally on 
$\cz$);

(c) the action (b) of $\cv$ on $\cz$ is transitive.
\endproclaim
Clearly, (a) is a consequence of 3.5(a).

We prove (b). We have $V^{N'}_{\ge a}=gV^N_{\ge a}=V^N_{\ge a}=V_{\ge a}$. Next we must
show that $\la gNg\i x,y\ra'+\la x,gNg\i y\ra'+\la gNg\i x,gNg\i y\ra'=0$ for 
$x,y\in V$ that is, $\la Ng\i x,g\i y\ra+\la g\i x,Ng\i y\ra+\la Ng\i x,Ng\i y\ra=0$ 
for $x,y\in V$. This follows from $\la Nx',y'\ra+\la x',Ny'\ra+\la Nx',Ny'\ra=0$ for 
$x',y'\in V$. Next we must show that $gNg\i,N$ induce the same map 
$\gr V_*@>>>\gr V_*$. (We must show: if $x\in V_{\ge a}$ then 
$gNg\i(x)-Nx\in V_{\ge a+3}$; this follows from $g\in\cv$.) Next we must show that for 
$x\in V_{\ge-a},y\in V_{\ge a}$ we have $\la x,y\ra'=\la x,y\ra$ that is 
$\la g\i x,g\i y\ra=\la x,y\ra$. Set $g\i=1+S$ where $S\in E_{\ge1}V_*$. We must show 
that $\la Sx,y\ra+\la x,Sy\ra+\la Sx,Sy\ra=0$. But $Sx\in V_{\ge1-a},y\in V_{\ge a}$ 
implies $\la Sx,y\ra=0$. Similarly $\la x,Sy\ra=0,\la Sx,Sy\ra=0$.

$\sp$ In the case where $p=2$ we see that the number $\nn$ defined in terms of $N,\lar$
is the same as that defined in terms of $N',\lar'$ and we must check that for 
$x\in V_{\ge-\nn}$ we have $\la x,(gNg\i)^{\nn-1}x\ra'=\la x,N^{\nn-1}x\ra$ 
that is, $\la g\i x,N^{\nn-1}g\i x\ra=\la x,N^{\nn-1}x\ra$ that is, 
$\la Sx,N^{\nn-1}x\ra+\la x,N^{\nn-1}Sx\ra+\la Sx,N^{\nn-1}Sx\ra=0$. We have 
$\la Sx,N^{\nn-1}x\ra+\la x,N^{\nn-1}Sx\ra=\la x,(N^{\nn-1}+(N^\da)^{\nn-1}Sx\ra$. This
is a linear combination of terms $\la x,N^{n'}Sx\ra$ where $n'\ge\nn$; each of these 
terms is zero since $x\in V_{\ge-\nn},N^{n'}Sx\in V_{\ge2n'-\nn+1}$ and 
$-\nn+2n'-\nn+1\ge1$. Next we have $\la Sx,N^{\nn-1}Sx\ra=0$ since 
$\la y,N^{\nn-1}y\ra=0$ for all $y\in V_{\ge1-\nn}$ by the definition of $\nn$. $\sp$

This completes the proof of (b). 

We show that (c) is a consequence of 3.6(c). Let 
$(N,\lar)\in\cz,(N',\lar')\in\cz$. By 2.6(a), since $V_*^N=V_*^{N'}$ and $N,N'$ induce 
the same $\nu$, we can find $S\in E_{\ge1}V_*$ such that $R=1+S$ satisfies $N'R=RN$. 
Define $\lar''\in\Symp(V)$ by $\la x,y\ra''=\la Rx,Ry\ra'$. From (b) we see that 
$(R\i N'R,\lar'')\in\cz$ that is $(N,\lar'')\in\cz$. Thus $\lar\in X_N,\lar''\in X_N$. 
By 3.6(c) we can find $S'\in E_{\ge1}V_*$ such that $R'=1+S'$ satisfies $R'N=NR'$ and 
$\la x,y\ra=\la R'x,R'y\ra''$ for all $x,y$ that is, $\la x,y\ra=\la RR'x,RR'y\ra'$. 
Then $RR'\in U'$ and $RR'N=RNR'=N'RR'$. Thus under the action (b), $RR'$ carries 
$(N,\lar)$ to $(N',\lar')$. This proves (c) (assuming 3.6(c)).

\subhead 3.8. Proof of 3.5(a)\endsubhead
We choose a direct sum decomposition $\op_{a\in\ZZ}V_a$ of $V$ as in 3.2(b). Define 
$N_2\in\End_2(V)$ by the requirement that $N_2:V_a@>>>V_{a+2}$ corresponds to
$\nu:\gr_aV_*@>>>\gr_{a+2}V_*$ under the obvious isomorphisms $V_a@>\si>>\gr_aV_*$,
$V_{a+2}@>\si>>\gr_{a+2}V_*$. 

$\sp$. If $p=2$ we regard $Q$ as a quadratic form on $V_{-\nn}$ via the obvious
isomorphism $V_{-\nn}@>\si>>\gr_{-\nn}V_*$. $\sp$

We will construct a linear map $N=\sum_{j\ge1}N_{2j}$ where $N_2$ is as above and for 
$j\ge2$, $N_{2j}\in\End(V)$ satisfy $N_{2j}V_a\sub V_{a+2j}$ for all $a$ and
$$\la\sum_{j\ge1}N_{2j}x,y\ra+\la x,\sum_{j\ge1}N_{2j}y\ra+
\la\sum_{j'\ge1}N_{2j'}x,\sum_{j''\ge1}N_{2j''}y\ra=0$$
for any $a,c$ and any $x\in V_a,y\in V_c$ that is,
$$\la N_{2j}x,y\ra+\la x,N_{2j}y\ra
+\sum_{j',j''\ge1;j'+j''=j}\la N_{2j'}x,N_{2j''}y\ra=0\tag a$$
for any $j\ge1$, any $a,c$ such that $a+c+2j=0$ and any $x\in V_a,y\in V_c$.

$\sp$ If $p=2$, we require in addition that $\la x,N^{\nn-1}x\ra=Q(x)$ for all
$x\in V_{-\nn}$ that is, $\sum_{i+i'=\nn-2}\la x,N_2^iN_4N_2^{i'}x\ra=Q(x)$ for all
$x\in V_{-\nn}$. $\sp$

We shall determine $N_j$ by induction on $j$. For $j=1$ the equation (a) is just 
$\la N_2x,y\ra+\la x,N_2y\ra=0$ for any $a,c$ such that $a+c+2=0$ and any 
$x\in V_a,y\in V_c$; this holds automatically by our choice of $N_2$. For $x\in V_a$ 
with $a<-2$ we set $N_4(x)=0$. Then the equation (a) for $j=2$ becomes

(b) $\la N_4x,y\ra+\la x,N_4y\ra=-\la N_2x,N_2y\ra$ for any $x\in V_{-2},y\in V_{-2}$,
$\la N_4x,y\ra=-\la N_2x,N_2y\ra$ for any $a>-2,x\in V_a,y\in V_{-a-4}$. 
\nl
The second equation in (b) determines uniquely $N_4(x)$ for $x\in V_a,a>-2$. Since
$\la N_2x,N_2y\ra$ is a symplectic form on $V_{-2}$ we can find $[,]\in\Bil(V_{-2})$ 
such that $[x,y]-[y,x]=-\la N_2x,N_2y\ra$ for any $x,y\in V_{-2}$. There is a unique 
linear map $N_4:V_{-2}@>>>V_2$ such that $\la N_4x,y\ra=[x,y]$ for any $x,y\in V_{-2}$.
Then equation (a) for $j=2$ is satisfied.

$\sp$ If $p=2$ the $N_4$ just determined satisfies
$\sum_{i+i'=\nn-2}\la x,N_2^iN_4N_2^{i'}x\ra=Q'(x)$ for all $x\in V_{-\nn}$, for some 
quadratic form $Q':V_{-\nn}@>>>\kk$ not necessarily equal to $Q$. For $x,y\in V_{-\nn}$
we have (by the choice of $N_4$):
$$\align&Q'(x+y)-Q'(x)-Q'(y)=\sum_{i+i'=\nn-2}\la x,N_2^iN_4N_2^{i'}y\ra+
\sum_{i+i'=\nn-2}\la y,N_2^iN_4N_2^{i'}x\ra\\&=\sum_{i+i'=\nn-2}\la N_2^ix,
N_4N_2^{i'}y\ra+\sum_{i+i'=\nn-2}\la N_4N_2^ix,N_2^{i'}y\ra
=\sum_{i+i'=\nn-2}\la N_2N_2^ix,N_2N_2^{i'}y\ra\\&=
\sum_{i+i'=\nn-2}\la x,N_2^\nn y\ra=\la x,N_2^\nn y\ra=Q(x+y)-Q(x)-Q(y).\endalign$$
It follows that $Q'(x)=Q(x)+\th(x)^2$ where $\th\in\Hom(V_{-\nn},\kk)$. We try to find 
$\z\in\End(V)$ with $\z(V_a)\sub V_{a+4}$ for all $a$ in such a way that (a) (for 
$j=2$) remains true when $N_4$ is replaced by $N_4+\z$ and
$\sum_{i+i'=\nn-2}\la x,N_2^i(N_4+\z)N_2^{i'}x\ra=Q(x)$ for $x\in V_{-\nn}$. (Then 
$N_4+\z$ will be our new $N_4$.) Thus we are seeking $\z$ such that

$\la\z(x),y\ra+\la x,\z(y)\ra=0$ for any $a,c$ with $a+c+4=0$ and $x\in V_a,y\in V_c$, 

$\sum_{i+i'=\nn-2}\la x,N_2^i\z N_2^{i'}x\ra=\th(x)^2$ for $x\in V_{-\nn}$.
\nl
The first of these two equations can be satisfied for $(a,c)\ne(-2,-2)$ by defining
$\z(x)=0$ for $x\in V_a,a\ne-2$. Then in the second equation the terms corresponding to
$i'$ such that $2i'-\nn\ne-2$ are $0$. Thus it remains to find a linear map
$\z:V_{-2}@>>>V_2$ such that

$\la \z(x),y\ra+\la x,\z(y)\ra=0$ for any $x,y\in V_{-2}$, 

$\la N_2^tx,\z N_2^tx\ra=\th(x)^2$ for $x\in V_{-\nn}$ where $t=(\nn-2)/2$.
\nl
Since $N_2^t:V_{-\nn}@>>>V_{-2}$ is injective (by the Lefschetz condition), there 
exists $\th_1\in\Hom(V_{-2},\kk)$ such that $\th_1(N_2^tx)=\th(x)$ for all 
$x\in V_{-\nn}$. We see that it suffices to find $\z\in\Hom(V_{-2},V_2)$ such that

$\la \z(x),y\ra+\la x,\z(y)\ra=0$ for any $x,y\in V_{-2}$, 

$\la x',\z x'\ra=\th_1(x')^2$ for $x'\in V_{-2}$.
\nl
It also suffices to find $b_0\in\Bil(V_{-2})$ such that $b_0=b_0^*$ and
$b_0(x,x)=\th_1(x)^2$ for $x\in V_{-2}$. Such $b_0$ clearly exists. $\sp$.

This completes the determination of $N_4$.

Now assume that $j\ge3$ and that $N_{2j'}$ is already determined for $j'<j$. For 
$x\in V_a$ with $a<-j$ we set $N_{2j}(x)=0$. Then equation (a) for our $j$ determines 
uniquely $N_{2j}(x)$ for $x\in V_a$ with $a>-j$. Next, we can find $[,]\in\Bil(V_{-j})$
such that

$[x,y]-[y,x]=-\sum_{j',j''\ge1;j'+j''=j}\la N_{2j'}x,N_{2j''}y\ra$.
\nl
To see this we observe that the right hand side is a symplectic form that is,
$\sum_{j',j''\ge1;j'+j''=j}\la N_{2j'}x,N_{2j''}x\ra=0$. There is a unique  
$N_{2j}\in\Hom(V_{-j},V_j)$ such that $\la N_{2j}x,y\ra=[x,y]$ for any $x,y\in V_{-j}$.
Then equation (a) for our $j$ is satisfied. This completes the inductive construction 
of $N$. We have $V_*^N=V_*$ by 2.4(a). We see that $N\in Y$. This completes the proof.

\subhead 3.9\endsubhead
In this subsection we prove 3.6(c) in a special case. Let $n\in\ZZ_{>0}$. We have 
$[-n,n]=I_0\sqc I_1$ where $I_\e=\{i\in[-n,n];i=\e\mod 2\}$ for $\e\in\{0,1\}$. For 
$i\in[-n,n]$ define $|i|\in\{0,1\}$ by $i=|i|\mod2$ that is by $i\in I_{|i|}$. Let 
$F_0,F_1\in\cc$. Let $V=\op_{i\in[-n,n]}F_i$ where $F_i=F_{|i|}$. A typical element of 
$V$ is of the form $(x_i)_{i\in[-n,n]}$ where $x_i\in F_{|i|}$. Define $N:V@>>>V$ by 
$(x_i)\m(x'_i)$ where $x'_i=x_{i-2}$ for $i\in[2-n,n]$, $x'_{-n}=0,x'_{1-n}=0$. We fix 
$\lar_0\in\Symp(V)$ such that 
$\la (x_i),(y_i)\ra_0=\sum_{i\in[-n,n]}(-1)^{(i-|i|)/2}b^{|i|}(x_i,y_{-i})$ where
$b^\e\in\Bil(F_\e)(\e\in\{0,1\}$ satisfy $b^{\e*}=(-1)^{1-\e}b^\e$, $b^\e$ is 
non-degenerate, $b^0\in\Symp(F_0)$. Note that $\la Nx,y\ra_0+\la x,Ny\ra_0=0$ for 
$x,y\in V$. 

We assume: if $p\ne2$ then either $F_0=0$ or $F_1=0$; $\sp$ if $p=2,b^1$ is symplectic
and $n\ge2$, then we are given a quadratic form $Q:F_0@>>>\kk$ such that 
$Q(x+y)=Q(x)+Q(y)+b^0(x,y)$ for $x,y\in F_0$. $\sp$

Let $X$ be the set of all $\lar\in\Symp(V)$ such that 
$\la Nx,y\ra+\la x,Ny\ra+\la Nx,Ny\ra=0$ for $x,y\in V$ and $\la x,y\ra=\la x,y\ra_0$ 
if there exists $i$ such that $x_j=0$ for $j\ne i$ and $y_j=0$ for $j\ne-i$; $\sp$ if 
$p=2,b^1$ is symplectic and $n\ge2$, we require also that $\la x,Nx\ra=Q(x_{-2})$ if 
$x\in V$ is such that $x_j=0$ for $j\ne-2$. $\sp$ 

Setting $\la (x_i),(y_i)\ra=\sum_{i,j}b_{ij}(x_i,y_j)$ identifies $X$ with the set of 
all families $(b_{ij})_{i,j\in[-n,n]}$ where $b_{ij}\in\Bil(F_{|i|},F_{|j|})$ are such 
that

$b_{i-2,j}+b_{i,j-2}+b_{ij}=0$ if $i,j\in[2-n,n]$,

$b_{i,-i}=(-1)^{(i-|i|)/2}b^{|i|}$ for all $i\in[-n,n]$,

$b_{ii}\in\Symp(F_{|i|})$ for all $i\in[-n,n]$, 

$b_{ij}^*=-b_{ji}$ for all $i,j\in[-n,n]$,

$b_{-2,0}(x,x)=Q(x)$ for $x\in F_0$ if $p=2,F_1=0$ and $n$ is even, $\ge2$.
\nl
(We have automatically $b_{ij}=0$ if $i+j\ge1$.)

Let $\D=\{T\in GL(V);TN=NT\}$, a subgroup of $GL(V)$; equivalently $\D$ is the set of 
linear maps $T:V@>>>V$ of the form 

(a) $T:(x_i)\m(x'_i),x'_i=\sum_{j\in[-n,i]}T_{i-j}^{|i|,|j|}x_j$
\nl
where $T_r^{\e,\d}\in\Hom(F_\d,F_\e)$ $(r\in[0,2n],\e,\d\in\{0,1\},r+\d=\e\mod 2)$ are 
such that $T_0^{00},T_0^{11}$ are invertible and $T_{2n}^{1-|n|,1-|n|}=0$. Now $\D$ 
acts on $X$ by $T:\lar\m\lar'$ where $\la Tx,Ty\ra'=\la x,y\ra$, or equivalently by 
$T:(b_{ij})\m(b'_{ij})$ where
$$b_{ij}(x,y)=\sum_{i'\in[i,n],j'\in[j,n]}b'_{i'j'}(T_{i'-i}^{|i'|,|i|}(x),
T_{j'-j}^{|j'|,|j|}(y)).$$
Let $\D_u=\{T\in\D;T_0^{00}=1,T_0^{11}=1\}$, a subgroup of $\D$. We show:

(b) {\it Let $k\in[1-n,0]$ and let $(\tb_{ij}),(b_{ij})$ be two points of $X$ such that
$b_{ij}=\tb_{ij}$ for $i+j\ge2k$. Then there exists $T\in\D_u$ such that 
$T(b_{ij})=(b'_{ij})$ and $b'_{ij}=\tb_{ij}$ for $i+j\ge2k-2$.}
\nl
For $\e\in\{0,1\}$ we set $a^\e=\tb_{ij}$ for $i,j\in[-n,n],i+j=-1,i=\e\mod2$. Then
$a^\e$ are independent of choices; they are $0$ unless $p=2$. We have 
$a^{\e*}=a^{1-\e}$. For $h\in\{2k-2,2k-1\}$ we set 
$c^\e_h=(-1)^{(i-\e)/2}(b_{ij}-\tb_{ij})$ where $i,j\in[-n,n],i+j=h,i=\e\mod2$. Then
$c^\e_h$ is independent of $i,j$. We have $c^\e_{2k-1}=0$ unless $p=2$. We have 
$c^{\e*}_{2k-2}=(-1)^{k-\e}c_{2k-2}^\e$, $c_{2k-1}^{\e*}=c_{2k-1}^{1-\e}$. Since 
$b_{k-1,k-1}-\tb_{k-1,k-1}$ is symplectic, $c_{2k-2}^\e$ is symplectic where 
$\e=k-1\mod2$.

{\it Case 1: $p\ne2$.} Let $\e\in\{0,1\}$ be such that $F_{1-\e}=0$. Since
$c^{\e*}_{2k-2}=(-1)^{k-\e}c_{2k-2}^\e$, we can find $\tc\in\Bil(F_\e)$ such that
$c_{2k-2}^\e=\tc+(-1)^{k-\e}\tc^*$. Since $b^\e$ is non-degenerate we can find
$\t\in\End(F_\e)$ such that $\tc(x,y)=b^\e(x,\t(y))$ for $x,y\in F_\e$. For 
$i,j\in[-n,n],i+j=2k-2,i=\e\mod2$ and $x,y\in F_\e$ we have
$$\align&b_{ij}(x,y)-\tb_{ij}(x,y)=(-1)^{(i-\e)/2}(\tc(x,y)+(-1)^{k-\e}\tc(y,x))\\&=
\tb_{i,j+2-2k}(x,\t(y))-\tb_{j,i+2-2k}(y,\t(x))=\tb_{i,j+2-2k}(x,\t(y))
+\tb_{i+2-2k,j}(\t(x),y).\endalign$$
Let $T$ be as in (a) with $T_0^{00}=1,T_0^{11}=1$, $T_{2-2k}^{\e,\e}=\t$ and the other 
components $0$. Define $(b'_{ij})$ by $T(b_{ij})=(b'_{ij})$. Then $(b'_{ij})$ has the 
required properties. 

$\sp$ {\it Case 2: $p=2,k=0$.} Since $b^0$ is non-degenerate we can find
$T_1^{0,1}\in\Hom(F_1,F_0)$ such that $c_{-1}^0(x,y)=\tb^0(x,T_1^{0,1}(y))$ for all 
$x\in F_0,y\in F_1$. Then $c_{-1}^1(x,y)=\tb^0(T_1^{0,1}(x),y)$ for all 
$x\in F_1,y\in F_0$. Thus, for $i\in I_0,j\in I_1,i+j=-1$ and $x\in F_0,y\in F_1$ we 
have $b_{ij}(x,y)+\tb_{ij}(x,y)=\tb_{i,-i}(x,T_1^{0,1}(y))$; for 
$i\in I_1,j\in I_0,i+j=-1$ and $x\in F_1,y\in F_0$ we have
$$b_{ij}(x,y)+\tb_{ij}(x,y)=\tb_{-j,j}(T_1^{0,1}(x),y).$$
Since $c_{-2}^{0*}=c_{-2}^0,b^{1*}=b^1$, we have $c_{-2}^0(y,y)=\th(y)^2$ for 
$y\in F_0$, $b^1(x,x)=\th_1(x)^2$ for $x\in F_1$ where $\th\in\Hom(F_0,\kk)$, 
$\th_1\in\Hom(F_1,\kk)$. If $b^1$ is not symplectic, we have $\th_1\ne0$. Hence there 
exists $T_1^{1,0}\in\Hom(F_0,F_1)$ such that $\th(y)=\th_1(T_1^{1,0}(y))$ for all 
$y\in F_0$. Then $c_{-2}^0(y,y)+b^1(T_1^{1,0}(y),T_1^{1,0}(y))=0$ for all $y\in F_0$. 
Thus $c_{-2}^0+b^1(T_1^{1,0}\ot T_1^{1,0})$ is symplectic. This also holds if $b^1$ is
symplectic (in that case we have 
$c_{-2}^0(y,y)=b_{-2,0}(y,y)-\tb_{-2,0}(y,y)=Q(y)-Q(y)=0$ for $y\in F_0$) and we take 
$T_1^{1,0}=0$. Now $c_{-2}^1$ is also symplectic.

Since $a^{0*}=a^1$, $a^1(T_1^{1,0}\ot1)+a^0(1\ot T_1^{1,0})$ is symplectic. Similarly 
$a^0(T_1^{0,1}\ot1)+a^1(1\ot T_1^{0,1})$ is symplectic. Hence 
$c_{-2}^0+b^1(T_1^{1,0}\ot T_1^{1,0})+a^1(T_1^{1,0}\ot1)+a^0(1\ot T_1^{1,0})$ is
symplectic and $c_{-2}^1+a^0(T_1^{0,1}\ot1)+a^1(1\ot T_1^{0,1})$ is symplectic. Hence 
we can find $\tc^0\in\Bil(F_0),\tc^1\in\Bil(F_1)$ such that 
$$c_{-2}^0+b^1(T_1^{1,0}\ot T_1^{1,0})+a^1(T_1^{1,0}\ot1)+a^0(1\ot T_1^{1,0})=
\tc^0+\tc^{0*},$$
$$c_{-2}^1+a^0(T_1^{0,1}\ot1)+a^1(1\ot T_1^{0,1})=\tc^1+\tc^{1*}.$$
Since $b^0,b^1$ are non-degenerate we can find $T_2^{0,0}\in\End(F_0)$,
$T_2^{1,1}\in\End(F_1)$ such that $\tc^0(x,y)=b^0(x,T_2^{0,0}(y))$ for $x,y\in F_0$, 
$\tc^1(x,y)=b^1(x,T_2^{1,1}(y))$ for $x,y\in F_1$. For $x,y\in F_0$ we have
$$\align&c_{-2}^0(x,y)+b^1(T_1^{1,0}(x)\ot T_1^{1,0}(x))+a^1(T_1^{1,0}(x),y)
+a^0(x,T_1^{1,0}(y))\\&=b^0(x,T_2^{0,0}(y))+b^0(T_2^{0,0}(x),y).\endalign$$
For $x,y\in F_1$ we have 
$$c_{-2}^1(x,y)+a^0(T_1^{0,1}(x),y)+a^1(x,T_1^{0,1}(y))=b^1(x,T_2^{1,1}(y))
+b^1(T_2^{1,1}(x),y).$$
Thus, for $i,j\in I_0,i+j=-2$ and $x,y\in F_0$ we have
$$\align&b_{ij}(x,y)=\tb_{ij}(x,y)+\tb_{i+1,j}(T_1^{1,0}(x),y)+
\tb_{i,j+1}(x,T_1^{1,0}(y))\\&+\tb_{i+1,j+1}(T_1^{1,0}(x),T_1^{1,0}(y))+
\tb_{i,-i}(x,T_2^{0,0}(y))+\tb_{-j,j}(T_2^{0,0}(x),y);\endalign$$
for $i,j\in I_1,i+j=-2$ and $x,y\in F_1$ we have
$$\align&b_{ij}(x,y)=\tb_{ij}(x,y)+\tb_{i+1,j}(T_1^{0,1}(x),y)+\tb_{i,j+1}(x,T_1^{0,1}
(y))\\&+\tb_{i,-i}(x,T_2^{1,1}(y))+\tb_{-j,j}(T_2^{1,1}(x),y).\endalign$$
Let $T$ be as in (a) with $T_0^{00}=1,T_0^{11}=1,T_1^{1,0},T_1^{0,1},T_2^{1,1}$,
$T_2^{0,0}$ as above and the other components $0$. Define 
$(b'_{ij})$ by $T(b_{ij})=(b'_{ij})$. Then $(b'_{ij})$ has the required properties. 

{\it Case 3: $p=2,k=-1$.} In this case we have $n\ge2$. We first show that there 
exists $\s\in\End(F_1)$ such that 
$$b^1(x,\s(y))=b^1(\s(x),y),\quad c^1_{-4}(x,x)=b^1(x,\s(x))+b^1(\s(x),\s(x))\tag *$$
for $x,y\in F_1$. The functions $F_1@>>>\kk,x\m b^1(x,x),x\m c^1_{-4}(x,x)$ are 
additive and homogeneous of degree $2$, hence are of the form 
$x\m\th(x)^2,x\m\th_1(x)^2$ where $\th,\th_1\in\Hom(F_1,\kk)$. We can find a direct sum
decomposition $F_1=F'\op F''$ where $b^1(x',x'')=0$ for all $x'\in F',x''\in F''$, 
$\th|_{F'}=0$, $F'=F_1$ if $\th=0$, $\dim F''\in\{1,2\}$ if $\th\ne0$. Define 
$\s'\in\End(F')$ by $\th_1(x)\th_1(y)=b^1(x,\s'(y))$ for $x,y\in F'$. Then 
$b^1(x,\s'(y))=b^1(\s'(x),y)$ for $x,y\in F'$, $\th_1(x)^2=b^1(x,\s'(x))+\th(\s'(x))^2$
for $x\in F'$.

If $\dim F''=1$ we have $\th|_{F''}\ne0$ and there is a unique $v\in F''$ such that 
$\th(v)=1$. Let $\s'':F''@>>>F''$ be multiplication by $a$ where $a\in\kk$ satisfies 
$a^2+a=\th_1(v)^2$. Then $\th_1(x)^2=b^1(x,\s''(x))+\th(\s''(x))^2$ for $x\in F''$ and 
$b(x,\s''(y))=b(\s''(x),y)$ for $x,y\in F''$.

If $\dim F''=2$ we can find a basis $\{v,v'\}$ of $F''$ such that 
$\th(v')=0,\th(v'')=1$. We set $b(v',v'')=f$. We have $f\ne0$. Define 
$\s''\in\End(F'')$ by $\s''(v')=\ta f\i v'+\ta v''$, $\s''(v'')=\th_1(v'')^2f\i v'$
where $\ta\in\kk$ satisfies $\ta^2+\ta=\th_1(v')^2$. Then $b(x,\s''(y))=b(\s''(x),y)$ 
for $x,y\in F''$, $\th_1(x)^2=b(x,\s''(x))+\th(\s''(x))^2$ for $x\in F''$. 

If $F''=0$ let $\s'':F''@>>>F''$ be the $0$ map. 

Define $\s\in\End(F_1)$ by $\s(x)=\s'(x)$ if $x\in F'$, $\s(x)=\s''(x)$ if $x\in F''$. 
Then $\s$ satisfies $(*)$. Since $b^0$ is non-degenerate we can find 
$T_3^{0,1}\in\Hom(F_1,F_0)$ such that
$$c^0_{-3}(x,y)+a^0(x,\s(y))=b^0(x,T_3^{0,1}(y))$$
for $x\in F_0,y\in F_1$. For any $i\in I_0,j\in I_1,i+j=-3$ and $x\in F_{|i|}$,
$y\in F_{|j|}$ we have
$$b_{ij}(x,y)=\tb_{ij}(x,y)+\tb_{i,j+2}(x,\s(y))+\tb_{i,j+3}(x,T_3^{0,1}(y)).$$
It follows that for any $i\in I_1,j\in I_0,i+j=-3$ and $x\in F_{|i|},y\in F_{|j|}$ we 
have
$$b_{ij}(x,y)=\tb_{ij}(x,y)+\tb_{i+2,j}(\s(x),y)+\tb_{i+3,j}(T_3^{0,1}(x),y).$$
Define $d_1\in\Bil(F_1)$ by $d_1(x,y)=\tb_{i,j+2}(x,\s(y))+\tb_{i+2,j}(\s(x),y)$ where 
$i,j\in I_1,i+j=-4$. Using the first equality in $(*)$ we see that $d_1$ is independent
of the choice of $i,j$. Define $d\in\Bil(F_1)$ by
$$d(x,y)=c^1_{-4}(x,y)+d_1(x,y)+b^1(\s(x),\s(y))+a^0(T_3^{0,1}(x),y)
+a^1(x,T_3^{0,1}(y)).$$
We have $d(x,x)=0$ for $x\in F_1$. (We use $(*)$ and the identity
$\tb_{i,j+2}+\tb_{j,i+2}=b^1$ for $i,j\in I_1,i+j=-4$.) Thus, $d$ is symplectic hence 
we can find $d'\in\Bil(F_1)$ such that $d=d'+d'{}^*$. Since $b^1$ is non-degenerate we 
can find $T_4^{1,1}\in\End(F_1)$ such that $d'(x,y)=b^1(x,T_4^{1,1}(y))$ for 
$x,y\in F_1$. We have
$$d(x,y)=b^1(x,T_4^{1,1}(y))+b^1(y,T_4^{1,1}(x))=b^1(x,T_4^{1,1}(y))
+b^1(T_4^{1,1}(x),y).$$
Hence for $i,j\in I_1,i+j=-4$ and $x,y\in F_1$ we have
$$\align&b_{ij}(x,y)=\tb_{ij}(x,y)+\tb_{i,j+2}(x,\s(y))+\tb_{i+2,j}(\s(x),y)\\&+
\tb_{i+2,j+2}(\s(x),\s(y))+\tb_{i+3,j}(T_3^{0,1}(x),y)+\tb_{i,j+3}(x,T_3^{0,1}(y))\\&
+b_{i,j+4}(x,T_4^{1,1}(y))+b_{i+4,j}(T_4^{1,1}(x),y).\endalign$$
For $i,j\in I_1,i+j=-2$ and $x,y\in F_1$ we have
$$b_{ij}(x,y)=\tb_{ij}(x,y)+\tb_{i,j+2}(x,\s(y))+\tb_{i+2,j}(\s(x),y)=\tb_{ij}(x,y).$$
Define $f\in\Bil(F_0)$ by $f(x,y)=c^0_{-4}(x,y)$. Then $f$ is symplectic. (We use that 
$c_{-4}^0$ is symplectic.) Hence we can find $f'\in\Bil(F_0)$ such that $f=f'+f'{}^*$. 
Since $b^0$ is non-degenerate we can find $T_4^{0,0}\in\End(F_0)$ such that
$f'(x,y)=b^0(x,T_4^{0,0}(y))$ for $x,y\in F_0$. We have
$$f(x,y)=b^0(x,T_4^{0,0}(y))+b^0(y,T_4^{0,0}(x))=b^0(x,T_4^{0,0}(y))
+b^0(T_4^{0,0}(x),y)$$
hence for $i,j\in I_0,i+j=-4$ and $x,y\in F_0$ we have
$$b_{ij}(x,y)=\tb_{ij}(x,y)+b_{i,j+4}(x,T_4^{0,0}(y))+b_{i+4,j}(T_4^{0,0}(x),y).$$
For $i,j\in I_0,i+j=-2$ and $x,y\in F_0$ we have
$$b_{ij}(x,y)=\tb_{ij}(x,y).$$
Let $T$ be as in (a) with $T_0^{00}=1,T_0^{11}=1$, $T_3^{0,1},T_4^{1,1},T_4^{0,0}$,
$T_2^{1,1}=\s$ as above and the other components $0$. Define $(b'_{ij})$ by 
$T(b_{ij})=(b'_{ij})$. Then $(b'_{ij})$ has the required properties. 

{\it Case 4: $p=2,k<-1$.} In this case we have $n\ge3$. Define $\e,\d\in\{0,1\}$ by 
$\e=k-1\mod 2,\d=1-\e$. Since $b^\d$ is non-degenerate, we have 
$c_{2k-2}^\d(x,y)=b^\d(x,\s(y))$ for $x,y\in F_\d$ where $\s\in\End(F_\d)$ is well
defined. Since $b^{\d*}=b^\d,c_{2k-2}^{\d*}=c_{2k-2}^\d$, we have 
$b^\d(x,\s(y))=b^\d(\s(x),y)$. Since $b^\e$ is non-degenerate we can find
$T_{1-2k}^{\e,\d}\in\Hom(F_\d,F_\e)$ such that
$$c^\e_{2k-1}(x,y)+a^\e(x,\s(y))=b^\e(x,T_{1-2k}^{\e,\d}(y))$$
for $x\in F_\e,y\in F_\d$. For any $i\in I_\e,j\in I_\d,i+j=2k-1$ and $x\in F_\e$,
$y\in F_\d$ we have
$$b_{ij}(x,y)=\tb_{ij}(x,y)+\tb_{i,j-2k}(x,\s(y))
+\tb_{i,j+1-2k}(x,T_{1-2k}^{\e,\d}(y)).$$
It follows that for any $i\in I_\d,j\in I_\e,i+j=2k-1$ and $x\in F_\d,y\in F_\e$ we 
have
$$b_{ij}(x,y)=\tb_{ij}(x,y)+\tb_{i-2k,j}(\s(x),y)
+\tb_{i+1-2k,j}(T_{1-2k}^{\e,\d}(x),y).$$
Define $d_1\in\Bil(F_\d)$ by $d_1(x,y)=\tb_{i,j-2k}(x,\s(y))+\tb_{i-2k,j}(\s(x),y)$ 
where $i,j\in I_\d,i+j=2k-2$. Using $b^\d(1\ot\s)=b^\d(\s\ot 1)$ we see that $d_1$ is 
independent of the choice of $i,j$. Define $d\in\Bil(F_\d)$ by
$$d(x,y)=c^\d_{2k-2}(x,y)+d_1(x,y)+a^\e(T_{1-2k}^{\e,\d}(x),y)
+a^\d(x,T_{1-2k}^{\e,\d}(y)).$$
We have $d(x,x)=0$ for $x\in F_\d$. (This follows from the choice of $\s$ and the 
identity $\tb_{i,j-2k}+\tb_{j,i-2k}=b^\d$ for $i,j\in I_\d,i+j=2k-2$.) Thus, $d$ is 
symplectic hence we can find $d'\in\Bil(F_\d)$ such that $d=d'+d'{}^*$. Since $b^\d$ is
non-degenerate we can find $T_{2-2k}^{\d,\d}\in\End(F_\d)$ such that 
$d'(x,y)=b^\d(x,T_{2-2k}^{\d,\d}(y))$ for $x,y\in F_\d$. We have
$$d(x,y)=b^\d(x,T_{2-2k}^{\d,\d}(y))+b^\d(y,T_{2-2k}^{\d,\d}(x))=
b^\d(x,T_{2-2k}^{\d,\d}(y))+b^\d(T_{2-2k}^{\d,\d}(x),y).$$
Hence for $i,j\in I_\d,i+j=2k-2$ and $x,y\in F_\d$ we have
$$\align&b_{ij}(x,y)=\tb_{ij}(x,y)+\tb_{i,j-2k}(x,\s(y))+\tb_{i-2k,j}(\s(x),y)+
\tb_{i+1-2k,j}(T_{1-2k}^{\e,\d}(x),y)\\&+\tb_{i,j+1-2k}(x,T_{1-2k}^{\e,\d}(y))
+b_{i,j+2-2k}(x,T_{2-2k}^{\d,\d}(y))+b_{i+2-2k,j}(T_{2-2k}^{\d,\d}(x),y).\endalign$$
For $i,j\in I_\d,i+j=2k$ and $x,y\in F_\d$ we have
$$b_{ij}(x,y)=\tb_{ij}(x,y)+\tb_{i,j-2k}(x,\s(y))+\tb_{i-2k,j}(\s(x),y)=\tb_{ij}(x,y).
$$
Define $f\in\Bil(F_\e)$ by $f(x,y)=c^\e_{2k-2}(x,y)$. Then $f$ is symplectic. (We use 
that $c^\e_{2k-2}$ is symplectic.) Hence we can find $f'\in\Bil(F_\e)$ such that 
$f=f'+f'{}^*$. Since $b^\e$ is non-degenerate we can find 
$T_{2-2k}^{\e,\e}\in\End(F_\e)$ such that $f'(x,y)=b^\e(x,T_{2-2k}^{\e,\e}(y))$ for 
$x,y\in F_\e$. We have
$$f(x,y)=b^\e(x,T_{2-2k}^{\e,\e}(y))+b^\e(y,T_{2-2k}^{\e,\e}(x))
=b^\e(x,T_{2-2k}^{\e,\e}(y))+b^\e(T_{2-2k}^{\e,\e}(x),y)$$
hence for $i,j\in I_\e,i+j=2k-2$ and $x,y\in F_\e$ we have
$$b_{ij}(x,y)=\tb_{ij}(x,y)+b_{i,j+2-2k}(x,T_{2-2k}^{\e,\e}(y))
+b_{i+2-2k,j}(T_{2-2k}^{\d,\d}(x),y).$$
For $i,j\in I_\d,i+j=2k$ and $x,y\in F_\e$ we have $b_{ij}(x,y)=\tb_{ij}(x,y)$.
Let $T$ be as in (a) with $T_0^{00}=1,T_0^{11}=1$, $T_{1-2k}^{\e,\d},T_{2-2k}^{\e,\e},
T_{2-2k}^{\d,\d}$, $T_{-2k}^{\d,\d}=\s$ as above and the other components $0$. Define 
$(b'_{ij})$ by $T(b_{ij})=(b'_{ij})$. Then $(b'_{ij})$ has the required properties. 
$\sp$ 

This completes the proof of (b).

We now verify the following special case of 3.6(c).

(c) {\it Let $(\tb_{ij}),(b_{ij})$ be two points of $X$. Then there exists $T\in\D_u$ 
such that $T(b_{ij})=(\tb_{ij})$.}
\nl
We first prove the following statement by induction on $k\in[-n,0]$.

($P_k$) {\it Assume in addition that $b_{ij}=\tb_{ij}$ for any $i,j$ with $i+j\ge2k$. 
Then there exists $T\in\D_u$ such that $T(b_{ij})=(\tb_{ij})$.}
\nl
If $k=-n$ the result is obvious. Assume now that $k\in[1-n,0]$. By (b) we can find 
$T'\in\D_u$ such that $T'(b_{ij})=(b'_{ij})$ and $b'_{ij}=\tb_{ij}$ for $i+j\ge2k-2$. 
By the induction hypothesis we can find $T''\in\D_u$ such that 
$T''(b'_{ij})=(\tb_{ij})$. Let $T=T''T'\in\D_u$. Then $T(b_{ij})=(\tb_{ij})$. This 
completes the proof of $(P_k)$ for $k\in[-n,0]$. In particular $(P_0)$ holds and (c) is
proved.

\subhead 3.10. Proof of 3.6(c)\endsubhead
Let $\lar,\lar'$ be two elements of the set $X$ in 3.6. We must show that $\lar,\lar'$ 
are in the same $U$-orbit. We argue by induction on $e$, the smallest integer $\ge0$ 
such that $N^e=0$. If $e=0$ we have $V=0$ and the result is obvious. If $e=1$ we have 
$N=0$. Then $V=\gr V_*$ canonically, $U=\{1\}$ and both $\lar,\lar'$ are the same as 
$\lar_0$ hence the result is clear. We now assume that $e\ge2$. 

$\sp$. Assume first that $p=2$. For $n\in\cl$ let $q_n:P^\nu_{-n}@>>>\kk$ be the 
quadratic forms attached to $(N,\lar)$ in 3.3 and let $q'_n:P^\nu_{-n}@>>>\kk$ be the 
analogous quadratic forms defined in terms of $(N,\lar')$. We show:

(a) {\it there exists $T\in U$ such that if $\lar''\in\Symp(V)$ is given by 
$\la x,y\ra''=\la Tx,Ty\ra$ then for $n\in\cl$ the quadratic form $q''_n$ defined as in
3.3 in terms of $(N,\lar'')$ satisfies $q''_n=q'_n$.}
\nl
We are seeking an $S\in E_{\ge1}V_*$ such that $SN=NS$ and
$\la(1+S)\dx,(1+S)N^{n-1}\dx\ra=\la\dx,N^{n-1}\dx\ra'$ that is, 
$\la(1+S)\dx,N^{n-1}(1+S)\dx\ra=\la\dx,N^{n-1}\dx\ra'$ that is, 

$\la S\dx,N^{n-1}\dx\ra+\la \dx,N^{n-1}S\dx\ra+\la S\dx,N^{n-1}S\dx\ra
=\la\dx,N^{n-1}\dx\ra'+\la\dx,N^{n-1}\dx\ra$
\nl
for any $n\in\cl$ and any $\dx\in V_{\ge-n}$ such that $N^{n+1}\dx=0$. Now 

$\la S\dx,N^{n-1}\dx\ra+\la\dx,N^{n-1}S\dx\ra=\la S\dx,(N^{n-1}+(N^\da)^{n-1})\dx\ra$
\nl
is a linear combination of terms $\la S\dx,N^{n'}\dx\ra$ with $n'\ge n$; each of these 
terms is $0$ since $S\dx\in V_{\ge1-n},N^{n'}\dx\in V_{\ge2n'-n}$ and $1-n+2n'-n\ge1$. 
Moreover, $\la S\dx,N^{n-1}S\dx\ra=\la \bS x,\nu^{n-1}\bS x\ra_0$ where 
$x\in P^\nu_{-n}$ is the image of $\dx$ and 
$\la\dx,N^{n-1}\dx\ra'+\la\dx,N^{n-1}\dx\ra=q'_n(x)+q_n(x)$. By the surjectivity of the
map $S\m\bS$ in 2.5(d), we see that it suffices to show that there exists 
$\s\in\End_1^\nu(\gr V_*)$ (that is $\s\in\End_1(\gr V_*)$ such that $\s\nu=\nu\s$) 
with $\la\s x,\nu^{n-1}\s x\ra_0=q'_n(x)+q_n(x)$ for any $n\in\cl$ and any 
$x\in P^\nu_{-n}$. 

For $n\in\cl'$ the last equation is automatically satisfied for any $\s$. (The left
hand side is zero by the definition of $\cl'$. The right hand side is equal by 3.3(c) 
to $Q'_\nn(\nu^{(n-\nn)/2}x)+Q_\nn(\nu^{(n-\nn)/2}x)$ where $Q_\nn$ is the quadratic 
form attached as in 3.3 to $(N,\lar)$ and $Q'_\nn$ is the analogous quadratic form 
defined in terms of $(N,\lar')$. The last sum is zero since $Q_\nn=Q'_\nn=Q$.)

We see that it suffices to show that there exists $\s\in\End_1^\nu(\gr V_*)$ such that
$\la\s x,\nu^{n-1}\s x\ra_0=q'_n(x)+q_n(x)$ for any $n\in\cl-\cl'$ and any 
$x\in P^\nu_{-n}$. 

For $n\in\cl-\cl'$, the quadratic forms $q'_n,q_n$ have the same associated symplectic
form (see 3.3(a)); hence there exists $\th_n\in\Hom(P^\nu_{-n},\kk)$ such that
$q'_n(x)+q_n(x)=\th_n(x)^2$ for all $x\in P^\nu_{-n}$. Hence it suffices to show that
the linear map 

$\r:\End_1^\nu(\gr V_*)@>>>\op_{n\in\cl-\cl'}\Hom(P^\nu_{-n},\kk)$
\nl
given by $\s\m(\th_n)$ where $\th_n(x)=\sqrt{\la \s x,\nu^{n-1}\s x\ra_0}$ for 
$x\in P^\nu_{-n}$ is surjective. Let $\ce=\op_{n\ge0}\Hom(P^\nu_{-n},\gr_{1-n}V_*)$. We
have an isomorphism $\p:\End_1^\nu(\gr V_*)@>\si>>\ce$ given by $\s\m(\s_n)$ where
$\s_n\in\Hom(P^\nu_{-n},\gr_{1-n}V_*)$ is the restriction of $\s$. Define a linear map

$\r':\ce@>>>\op_{n\in\cl-\cl'}\Hom(P^\nu_{-n},\kk)$
\nl
by $(\s_n)\m(\th_n)$ where $\th_n(x)=\sqrt{\la \s_nx,\nu^{n-1}\s_nx\ra_0}$ for 
$x\in P^\nu_{-n}$. We have $\r'\p=\r$. Hence it suffices to show that $\r'$ is 
surjective. It also suffices to show that for any $n\in\cl-\cl'$ the linear map

$\r'_n:\Hom(P^\nu_{-n},\gr_{1-n}V_*)@>>>\Hom(P^\nu_{-n},\kk)$
\nl
given by $f\m\th$, where $\th(x)=\sqrt{\la fx,\nu^{n-1}fx\ra_0}$ for $x\in P^\nu_{-n}$,
is surjective. Define $g\in\Hom(\gr_{1-n}V_*@>>>\kk)$ by 
$h\m\sqrt{\la h,\nu^{n-1}h\ra_0}$. Then $\r'_n(f)=g\circ f$ for 
$f\in\Hom(P^\nu_{-n},\gr_{1-n}V_*)$. Hence to show that $\r'_n$ is surjective it 
suffices to show that $g\ne0$. Since $n\in\cl-\cl'$, there exists $m$ odd such that 
$m\ge n+3$ and $b_m$ is not symplectic. Hence there exists $u'\in P^\nu_{-m}$ such that
$\la u',\nu^mu'\ra_0\ne0$. We have $m=(n-1)+2k$ where $k$ is an integer $\ge2$. Let 
$u=N^ku'\in\gr_{1-n}V_*$ and 

$\la u,\nu^{n-1}u\ra_0=\la \nu^ku',\nu^{n-1+k}u'\ra_0=\la u',\nu^mu'\ra_0\ne0$.
\nl
Thus $g(u)\ne0$. We see that $g\ne0$, as required. This proves (a).

Note that $\lar''$ in (a) is in $X$ (in fact in the $U$-orbit of $\lar$). Replacing if 
necessary $\lar$ by $\lar''$ we see that 

(b) {\it we may assume that $\lar,\lar'$ are such that $q_n=q'_n$ for all $n\in\cl$.} 
$\sp$

We now return to a general $p$. Let $r\ge e$. Let $F$ be a complement of 
$V_{\ge2-r}=\ker N^{r-1}$ in $V_{\ge1-r}=V$ and let $F'$ be a complement of 
$V_{\ge3-r}=\ker N^{r-2}+NV$ in $V_{\ge2-r}=\ker N^{r-1}$. Consider the linear map $\a$
of $F\op F'\op F\op\do\op F'\op F$ ($2r-1$ summands) into $V$ given by

$(x_{1-r},x_{2-r},\do,x_{r-2},x_{r-1})\m$

$x_{1-r}+Nx_{3-r}+\do+N^{r-1}x_{r-1}+x_{2-r}+Nx_{4-r}+\do+N^{r-2}x_{r-2}$.
\nl
(Here $x_i\in F$ if $i=r+1\mod2$ and $x_i\in F'$ if $i=r\mod2$.) Let $W$ be the image 
of $\a$. We show that 

(c) $\lar$ and $\lar'$ are non-degenerate on $W$.
\nl 
We prove this only for $\lar$; the proof for $\lar'$ is the same. Assume that 
$w=x_{1-r}+Nx_{3-r}+\do+N^{r-1}x_{r-1}+x_{2-r}+Nx_{4-r}+\do+N^{r-2}x_{r-2}$ with $x_i$ 
as above satisfies $\la w,W\ra=0$. We show that each $x_i$ is $0$. We have 
$0=\la w,N^{r-1}F\ra=\la x_{1-r},N^{r-1}F\ra=0$. Using the non-degeneracy of $b_{r-1}$ 
we see that $x_{1-r}=0$ and 
$w=Nx_{3-r}+\do+N^{r-1}x_{r-1}+x_{2-r}+Nx_{4-r}+\do+N^{r-2}x_{r-2}$. We have 
$0=\la w,N^{r-2}F'\ra=\la x_{2-r},N^{r-2}F'\ra$. Using the non-degeneracy of $b_{r-2}$ 
we see that $x_{2-r}=0$ and 
$w=Nx_{3-r}+\do+N^{r-1}x_{r-1}+Nx_{4-r}+\do+N^{r-2}x_{r-2}$. We have 
$0=\la w,N^{r-2}F\ra=\la Nx_{3-r},N^{r-2}F\ra=-\la x_{3-r},N^{r-1}F\ra$. Using the 
non-degeneracy of $b_{r-1}$ we see that $x_{3-r}=0$. Continuing in this way we see that
each $x_i$ is $0$. This proves (c).

The proof shows also that $\a$ is injective. 

Let $Z=\{x\in V;\la x,W\ra=0\},Z'=\{x\in V;\la x,W\ra'=0\}$. From (c) we see that
$V=W\op Z=W\op Z'$. 

Clearly, $W$ is $N$-stable hence $(1+N)$-stable. Since $1+N$ is an isometry of $\lar$ 
it follows that $Z$ is $(1+N)$-stable hence $N$-stable. Similarly, $Z'$ is $N$-stable. 
Define $\Ph\in GL(V)$ by $\Ph(x)=x$ for $x\in W$, $\Ph(x)=x'$ for $x\in Z$ where 
$x'\in Z'$ is given by $x-x'\in W$. We have $\Ph\in U$ (see 2.7(c),(d)). Define 
${}'\lar\in\Symp(V)$ by ${}'\la x,y\ra=\la \Ph(x),\Ph(y)\ra'$. By 3.6(a), we have 
${}'\lar\in X$.

Let ${}'Z=\{x\in V;{}'\la x,W\ra=0\}$. We show that ${}'Z=Z$. Let $x=x_1+x_2$ where 
$x_1\in W,x_2\in Z$. We have $x_2=w+x'_2$, $w\in W,x'_2\in Z'$. For $w'\in W$ we have 
$\la\Ph(x),w'\ra'=\la x_1+x'_2,w'\ra'=\la x_1,w'\ra'$. The condition that 
$\la\Ph(x),W\ra'=0$ is that $\la x_1,W\ra'=0$ or that $x_1=0$ (using (c)) or that 
$x\in Z$. Thus, 
${}'Z=\{x\in V;\la \Ph(x),\Ph(W)\ra'=0\}=\{x\in V;\la \Ph(x),W\ra'=0\}=Z$ as required.

$\sp$ In the case where $p=2$ we show that for any $n\in\cl$ the quadratic form $q'_n$
attached to $(N,\lar')$ as in 3.3 is equal to the analogous quadratic form attached to 
$(N,{}'\lar)$. We must show that, if $x\in V_{\ge-n}$, $N^{n+1}x=0$ then
$\la\Ph x,\Ph N^{n-1}x\ra'=\la x,N^{n-1}x\ra'$ that is, 
$\la\Ph x,N^{n-1}\Ph x\ra'=\la x,N^{n-1}x\ra'$. Both sides are additive in $x$. We can 
write $x=x_1+x_2$ where $x_1\in W,x_2\in Z$ satisfy $x_1,x_2\in V_{\ge-n}$,
$N^{n+1}x_1=0,N^{n+1}x_2=0$. We may assume that $x=x_1$ or $x=x_2$. When $x=x_1$ the 
desired equality is obvious. Hence we  may assume that $x\in Z$. Write $x=x'+w$,
$x'\in Z',w\in W$. We must show that $\la x+w,N^{n-1}x+N^{n-1}w\ra'=\la x,N^{n-1}x\ra'$
that is, $\la x,N^{n-1}w\ra'+\la w,N^{n-1}x\ra'+\la w,N^{n-1}w\ra'=0$ that is,
$\la x,(N^{n-1}+(N^\da)^{n-1})w\ra'+\la w,N^{n-1}w\ra'=0$ that is, 
$\la x,N^nw\ra'+\la w,N^{n-1}w\ra'=0$ (we use $N^{n+1}w=0$) that is, 
$\la x'+w,N^nw\ra'+\la w,N^{n-1}w\ra'=0$ that is, 
$\la w,N^nw\ra'+\la w,N^{n-1}w\ra'=0$. Now $w\in W_{\ge1-n}$ (see 2.7(b)), 
$N^nw\in W_{\ge n+1}$ hence $\la w,N^nw\ra'=0$. It remains to show 
$\la w,N^{n-1}w\ra'=0$. Since $w\in W_{\ge1-n}$, $N^{n+1}w=0$, it suffices to show 
$\la y,\nu^{n-1}y\ra_0=0$ for any $y\in\gr_{1-n}V_*$ such that $\nu^{n+1}y=0$. This has
already been seen in the proof in 3.3 that $q_n$ is well defined. $\sp$

Replacing $\lar'$ by ${}'\lar$ (which is in the same $U$-orbit) we see that condition 
(b) is preserved (for $p=2$).

Thus, we may assume that $\lar,\lar'$ satisfy $Z=Z'$ and that for $p=2$ condition (b) 
holds. Thus $V=W\op Z$ is an othogonal decomposition with respect to either $\lar$ or 
$\lar'$. Let $\lar_W,\lar_Z$ be the restrictions of $\lar$ to $W,Z$. Let 
$\lar'_W,\lar'_Z$ be the restrictions of $\lar'$ to $W,Z$. Let $U_1$ (resp. $U_2$) be 
the analogue of $U$ for $W$ (resp. $Z$) defined in terms of $N$ and $W^N_*$ (resp. 
$Z^N_*$). We have naturally $U_1\T U_2\sub U$.

We consider $5$ cases.

{\it Case 1: $p\ne2$.} Take $r=e+1$. (Thus, $F=0$.) By the induction hypothesis, 
$\lar_Z$ is carried to $\lar'_Z$ by some $u_2\in U_2$. By 3.9, $\lar_W$ is carried to 
$\lar'_W$ by some $u_1\in U_1$. Then $\lar$ is carried to $\lar'$ by $(u_1,u_2)\in U$. 

$\sp$ {\it Case 2: $p=2$, $e$ is odd and $b_{e-2}$ is symplectic.} Take $r=e+1$. (Thus,
$F=0$.) We have $e-1\in\cl$. The sets $\cl$ attached to $\lar_Z,\lar'_Z$ are the same 
as $\cl$ for $\lar,\lar'$. The quadratic forms attached to $\lar_Z$, $\lar'_Z$ and 
$n\in\cl-\{e-1\}$ are the same as those attached to $\lar,\lar'$ and $n$, hence they 
coincide. The quadratic forms attached to $\lar_Z,\lar'_Z$ and $n=e-1$ also coincide: 
they are both $0$. Hence the Quadratic forms attached to $\lar_Z,\lar'_Z$ coincide (see
3.3(c)). The quadratic forms attached to $\lar_W,\lar'_W$ coincide: for $e-1$ they are 
the same as those attached to $\lar,\lar'$ and $e-1$ and for other $n$ they are zero. 
Hence the Quadratic forms attached to $\lar_W,\lar'_W$ coincide. By the induction 
hypothesis, $\lar_Z$ is carried to $\lar'_Z$ by some $u_2\in U_2$. By 3.9, $\lar_W$ is 
carried to $\lar'_W$ by some $u_1\in U_1$. Then $\lar$ is carried to $\lar'$ by 
$(u_1,u_2)\in U$. 

{\it Case 3: $p=2$, $e$ is even and $b_{e-1}$ is symplectic.} Take $r=e+1$. (Thus, 
$F=0$.) The sets $\cl$ attached to $\lar_Z,\lar'_Z$ are the same as $\cl$ for 
$\lar,\lar'$. The quadratic forms attached to $\lar_Z,\lar'_Z$ and $n\in\cl$ are the 
same as those attached to $\lar,\lar'$ and $n$, hence they coincide. Hence the 
Quadratic forms attached to $\lar_Z,\lar'_Z$ coincide. By the induction hypothesis, 
$\lar_Z$ is carried to $\lar'_Z$ by some $u_2\in U_2$. By 3.9, $\lar_W$ is carried to 
$\lar'_W$ by some $u_1\in U_1$. Then $\lar$ is carried to $\lar'$ by $(u_1,u_2)\in U$. 

{\it Case 4: $p=2$, $e$ is even and $b_{e-1}$ is not symplectic.} Take $r=e$. The sets 
$\cl$ attached to $\lar_Z,\lar'_Z$ are the same as $\cl$ for $\lar,\lar'$. The 
quadratic forms attached to $\lar_Z,\lar'_Z$ and $n\in\cl$ are the same as those 
attached to $\lar,\lar'$ and $n$, hence they coincide. Hence the Quadratic forms 
attached to $\lar_Z,\lar'_Z$ coincide. By the induction hypothesis, $\lar_Z$ is carried
to $\lar'_Z$ by some $u_2\in U_2$. By 3.9, $\lar_W$ is carried to $\lar'_W$ by some 
$u_1\in U_1$. Then $\lar$ is carried to $\lar'$ by $(u_1,u_2)\in U$. 

{\it Case 5: $p=2$, $e$ is odd, $\ge3$ and $b_{e-2}$ is not symplectic.} Take $r=e$. By
3.9, $\lar_W$ is carried to $\lar'_W$ by some $u_1\in U_1$. Replacing $\lar'$ by a 
translate under $(u_1,1)\in U$ we see that we may assume in addition that 
$\lar_W=\lar'_W$. Let $\tW=F+NF+\do+N^{r-1}F$. Let 
$W'=\{w\in W;\la w,\tW\ra=0\}=\{w\in W;\la w,\tW\ra'=0\}$. Then $W=\tW\op W'$, 
orthogonal direct sum for both $\lar,\lar'$. Let $\tZ=W'\op Z$, orthogonal direct sum 
for both $\lar,\lar'$. Then $V=\tW\op\tZ$, orthogonal direct sum for both $\lar,\lar'$.
Let $\lar_{\tZ},\lar'_{\tZ}$ be the restrictions of $\lar,\lar'$ to $\tZ$. Let $\tU_1$ 
(resp. $\tU_2$) be the analogue of $U$ for $\tW$ (resp. $\tZ$) defined in terms of $N$ 
and $\tW^N_*$ (resp. $\tZ^N_*$). We have naturally $\tU_1\T\tU_2\sub U$. The sets $\cl$
attached to $\lar_{\tZ}$, $\lar'_{\tZ}$ are the same as $\cl$ for $\lar,\lar'$. The 
quadratic forms attached to $\lar_{\tZ},\lar'_{\tZ}$ and $n\in\cl$ are the same as 
those attached to $\lar,\lar'$ and $n$, hence they coincide. Hence the Quadratic forms 
attached to $\lar_{\tZ},\lar'_{\tZ}$ coincide. By the induction hypothesis, 
$\lar_{\tZ}$ is carried to $\lar'_{\tZ}$ by some $\tu_2\in\tU_2$. Then $\lar$ is 
carried to $\lar'$ by $(1,\tu_2)\in U$. $\sp$

This completes the proof of 3.6(c) hence also that of 3.5, 3.6, 3.7.

\subhead 3.11\endsubhead
Here is the order of the proof of the various assertions in Propositions 3.5-3.7:
3.5(a) (see 3.8); 3.7(a) (see 3.7); 3.6(a) (see 3.6); 3.7(b) (see 3.7); 3.6(b) (see 
3.6); 3.5(b) (see 3.5); 3.6(c) (see 3.9, 3.10); 3.7(c) (see 3.7); 3.5(c) (see 3.5).

\subhead 3.12\endsubhead
Let $V\in\cc$ and let $\lar\in\Symp(V)$. The following result can be deduced from 
\cite{\SPA, I, 2.10}.

Let $C,C_0$ be two $GL(V)$-conjugacy classes in $\Nil(V)$ such that 
$C\cap\cm_{\lar}\ne\em,C_0\cap\cm_{\lar}\ne\em$ and $C$ is contained in the closure of 
$C_0$ in $GL(V)$. Then $C\cap\cm_{\lar}$ is contained in te closure of 
$C_0\cap\cm_{\lar}$ in $\cm_{\lar}$.

\subhead 3.13\endsubhead
Let $V\in\cc$ and let $\lar\in\Symp(V)$. Let $G=Sp(\lar)$. For any self-dual filtration
$V_*$ of $V$ and for $n\ge1$ let $E^{\lar}_{\ge n}V_*=E_{\ge n}V_*\cap\cm_{\lar}$, a 
unipotent algebraic group with multiplication $T*T'=T+T'+TT'$. Let
$$\tix(V_*)=\x(V_*)\cap\cm_{\lar}=\{N\in\cm_{\lar};V^N_*=V_*\}
=\{N\in E^{\lar}_{\ge2}V_*;\bN\in\End_2^0(\gr V_*)\}$$
(see 2.9). The following three conditions are equivalent:

(i) $\tix(V_*)\ne\em$;

(ii) there exists $\nu\in\End_2^0(\gr V_*)$ which is skew-adjoint with respect to the 
symplectic form on $\gr V_*$ induced by $\lar$;

(iii) $\dim\gr_nV_*=\dim\gr_{-n}V_*\ge\dim\gr_{-n-2}V_*$ for all $n\ge0$ and
$\dim\gr_{-n}V_*=\dim\gr_{-n-2}V_* \mod2$ for all $n\ge0$ even.
\nl
We have (i)$\imp$(ii) by the definition of $\tix(V_*)$; we have (ii)$\imp$(iii) by
2.3(d) and 3.1(c). Now (iii)$\imp$(ii) is easily checked. We have (ii)$\imp$(iii) by 
3.5(a).

Let $\fF_{\lar}$ be the set of all self-dual filtrations $V_*$ of $V$ that satisfy 
(i)-(iii). From the definitions we have a bijection 

(a) $\fF_{\lar}@>\si>>D_G,V_*\m\l$
\nl
($D_G$ as in 1.1) where $\l=(G^\l_0\sps G^\l_1\sps G^\l_2\sps\do)$ is defined in terms 
of $V_*$ by $G^\l_0=E_{\ge0}V_*\cap G$ and $G^\l_n=1+E_{\ge n}^{\lar}V_*$ for $n\ge1$. 
The sets $\tix(V_*)$ (with $V_*\in\fF_{\lar})$ form a partition of $\cm_{\lar}$. (If 
$N\in\cm_{\lar}$ we have $N\in\tix(V_*)$ where $V_*=V^N_*$.)

Let $V_*\in\fF_{\lar}$. Let $C_0$ be the unique $GL(V)$-conjugacy class in $\Nil(V)$
that contains $\x(V_*)$. We have 

$E_{\ge2}^{\lar}V_*-\tix(V_*)=(E_{\ge2}V_*-\x(V_*))\cap\cm_{\lar}
=E_{\ge2}V_*\cap(\cup_CC)\cap\cm_{\lar}$
\nl
(the last equality follows from 2.9; $C$ runs over all $GL(V)$-orbits in $\Nil(V)$ such
that $C\sub\bC_0-C_0$). Using 3.12 we see that

$E_{\ge2}^{\lar}V_*-\tix(V_*)=E_{\ge2}^{\lar}V_*\cap(\cup_C(C\cap\cm_{\lar}))$
\nl
where $C$ runs over all $GL(V)$-orbits in $\Nil(V)$ such that $C\cap\cm_{\lar}\ne\em$ 
and $C\sub\ov{C_0\cap\cm_{\lar}}-(C_0\cap\cm_{\lar})$. We see that, if $V_*\m\l$ (as in
(a)) and $\bla$ is the $G$-orbit of $\l$ in $D_G$ then (with notation of 1.1) 
$\tH^\bla$ is the union of $G$-conjugacy classes in $1+\cm_{\lar}$ contained in 
$1+\bC_0$, $H^\bla$ is the union of $G$-conjugacy classes in $1+\cm_{\lar}$ contained 
in $1+C_0$, $X^\l=1+\tix(V_*)=1+(E_{\ge2}^{\lar}V_*\cap C_0)$. We see that 
$\fP_1-\fP_3$ hold.

\subhead 3.14\endsubhead
We preserve the setup of 3.13. Let $V_*\in\fF_{\lar}$ and let $\l\in D_G$ be the 
corresponding element. Define $\lar_0\in\Symp(\gr V_*)$ as in 3.2. The map 
$E_{\ge2}^{\lar}V_*@>>>\End_2^0(\gr V_*),N\m\bN$ restricts to a map 

$\p:\tix(V_*)@>>>E:=\{\nu\in\End_2^0(\gr V_*);\nu\text{ skew-adjoint with respect to }
\lar_0\}$. 
\nl
We show:

(a) {\it The group $E_{\ge3}^{\lar}V_*$ (see 3.13) acts freely on $\tix(V^*)$ by
$T,N\m T*N$ (see 3.13) and the orbit space of this action may be identified with $E$
via $\p$.}
\nl
We show this only at the level of sets. If $T\in E_{\ge3}^{\lar}V_*,N\in\tix(V_*)$ then
$T*N\in E_{\ge2}^{\lar}V_*$ and $T*N,N$ induce the same map in $\End_2(\gr V_*)$; hence
$T*N\in\tix(V_*)$. Thus $T,N\m T*N$ is an action of $E_{\ge3}^{\lar}V_*$ on 
$\tix(V_*)$. This action is free: it is the restriction of the action of 
$E_{\ge3}^{\lar}V_*$ on $E_{\ge2}^{\lar}V_*$ by left multiplication for the group 
structure in 3.13. If $N,N'\in\tix(V_*)$ induce the same map in $\End_2^0(\gr V_*)$ 
then $N'-N\in E_{\ge3}V_*$. Set $T=(N'-N)(1+N)\i\in E_{\ge3}V_*$. Then 
$(1+T)(1+N)=1+N'$ and we have automatically $T\in E^{\lar}_{\ge3}V_*$ and $T+N=N'$. 
Thus the orbits of the $E_{\ge3}^{\lar}V_*$-action on $\ti\x(V^*)$ are exactly the 
non-empty fibres of $\p$. It remains to show that $\p$ is surjective. This follows from
3.5(a).

Now let $N,N'\in\tix(V_*)$ be such that $\bN=\bN'=\nu\in\End_2^0(\gr V_*)$. We show:

(b) {\it there exists $g\in E_{\ge0}V_*\cap G$ such that $N'=gNg\i$.}
\nl
Assume first that $p=2$. The set $\cl\sub2\NN$ defined in 3.3 in terms of $N$ is the 
same as that defined in terms of $N'$. Let $q_n:P^\nu_{-n}@>>>\kk$ be the quadratic
form defined in terms of $N$ (for $n\in\cl$) as in 3.3 and let $q'_n:P^\nu_{-n}@>>>\kk$
be the analogous quadratic form defined in terms of $N'$. From 3.3(a) we see that for 
any $n\in\cl$ there exists an automorphism $h_n:P^\nu_{-n}@>>>P^\nu_{-n}$ which 
preserves the symplectic form $x,y\m b_n(x,y)$ (see 3.1) and satisfies
$q'_n(x)=q_n(h_nx)$ for any $x\in P^\nu_{-n}$. There is a unique $h\in Sp(\lar_0)$ such
that $h(x)=h_n(x)$ for $x\in P^\nu_{-n},n\in\cl$, $h(x)=x$ for $x\in P^\nu_{-n}$,
$n\in\ZZ-\cl$, $h\nu=\nu h$. Let $V=\op_aV_a$ be a direct sum decomposition as in 
3.2(b). Then $\End_0(V)$ is defined and we define $\tih\in\End_0(V)$ by the requirement
that for any $a$, $\tih:V_a@>>>V_a$ corresponds to $h:\gr_aV_*@>>>\gr_aV_*$ under the 
obvious isomorphism $V_a@>\si>>\gr_aV_*$. Then $\tih\in E_{\ge0}V_*\cap G$ and 
$\tih N\tih\i=N''$ where $N''\in E^{\lar}_{\ge2}V_*$ satisfies $\bN''=\nu$. Moreover, 
the quadratic form $P^\nu_{-n}@>>>\kk$ defined as in 3.3 in terms of $N''$ (instead of 
$N$) for $n\in\cl$ is $x\m h_n(x)$ that is, $q'_n$. From 3.3(c) we see that the 
Quadratic form $Q_n$ defined for $n\in\cl'$ in terms of $N''$ is the same as that 
defined in terms of $N'$. From 3.5(c) we see that there exists 
$h'\in1+E^{\lar}_{\ge1}V_*$ such that $h'N''h'{}\i=N'$. Setting 
$g=h'\tih\in E_{\ge0}V_*\cap G$ we have $gNg\i=N'$.

Next assume that $p\ne2$. From 3.5(c) we see that there exists 
$g\in1+E^{\lar}_{\ge1}V_*$ such that $gNg\i=N'$. This proves (b).

We see that $\fP_6$ (hence $\fP_4$) holds.

From (a) we see that the $G^\l_0$-action on $\tix(V^*)$ (conjugation) induces an action
of $\bG^\l_0=G^\l_0/G^\l_1$ on $E$ and from (b) we see that this gives rise to a
bijection between the set of $G^\l_0$-orbits on $\tix(V^*)$ and the set of
$\bG^\l_0$-orbits on $E$. We describe this last set of orbits. We identify $\bG^\l_0$ 
with $\End_0(\gr V_*)\cap Sp(\lar_0)$ with the action on $E$ given by $g:\nu\m\nu'$ 
where $\nu'(x)=g\nu(g\i x)$ for $x\in\gr V_*$. 

Let $I=\{n\in2\NN+1,\dim\gr_{-n}V_*-\dim\gr_{-n-2}V_*\in\{2,4,6,\do\}\}$. For any 
subset $J\sub I$ let $E_J$ be the set of all $\nu\in E$ such that for any $n\in I$ we 
have

$\{x\in\gr_{-n}V_*;\nu^{n+1}x=0,\la x,\nu^nx\ra_0\ne0\}\ne\em\lra n\in J$.
\nl
Let $\bE$ be the set of all direct sum decompositions $\gr V_*=\op_{n\ge0}W^n$ where 
$W^n\in\bcc$ (see 2.1) are such that $\la W^n,W^{n'}\ra_0=0$ for $n\ne n'$ and for 
$n\ge 0$, $\dim W^n_a$ is $\dim\gr_{-n}V_*-\dim\gr_{-n-2}V_*$ if $a\in[-n,n],a=n\mod2$ 
and is $0$ for other $a$. Define $\ph:E@>>>\bE$ by $\nu\m(W^n)$ where 
$W^n=\sum_{k\ge0}\nu^kP^\nu_{-n}$. Then $\ph$ is $\bG^\l_0$-equivariant where 
$\bG^\l_0$ acts on $\bE$ in an obvious way (transitively).

Let $w=(W^n)\in\bE$. Let $G^w$ be the stabilizer of $w$ in $\bG^\l_0$. Let 
$E^w=\ph\i(w)$. Now $E^w$ may be identified with $\prod_{n\ge0}E^w_n$ where $E^w_n$ is 
the set of all skew-adjoint elements in $\End_2^0(W^n)$ with respect to 
$\lar_0|_{W^n}$. Moreover $G^w$ may be identified with $\prod_{n\ge0}G^w_n$ where 
$G^w_n=\End_0(W^n)\cap Sp(\lar_0|_{W^n})$. Furthermore, we may identify 
$E^w_n=E^{w1}_n\T E^{w2}_n$ where $E^{w1}_n$ consists of all sequences of isomorphisms 

(c) $W^n_{-n}@>\si>>W^n_{-n+2}@>\si>>W^n_{-n+4}@>\si>>\do@>\si>>W^n_{-\d}$
\nl
($\d=0$ if $n$ is even and $\d=1$ if $n$ is odd) and $E^{w2}_n$ is the set of 
non-degenerate symmetric bilinear forms $W^n_{-1}\T W^n_{-1}@>>>\kk$ (if $n$ is odd) 
and is a point if $n$ is even. (This identification is obtained by attaching to 
$\nu\in E^w_n$ the isomorphisms (c) induced by $\nu$ and if $n$ is odd, the bilinear 
form $x,x'\m\la x,\nu x'\ra_0$ on $W^n_{-1}$.)

We claim that if $p=2$, the subsets $E_J$ are precisely the orbits of $\bG^\l_0$ on $E$
while if $p\ne2$, $E$ is a single orbit of $\bG^\l_0$. Using the transitivity of the
$\bG^\l_0$ action on $\bE$ we see that it suffices to prove: if $p=2$, the subsets 
$E^w_J=E_J\cap E^w$ are precisely the $G^w$-orbits on $E$ while if $p\ne2$, $E^w$ is a
single $G^w$-orbit. If $n\n I$, $G'_n$ acts transitively on $E'_n$. If $n\in I$, 
$pr_2:E^w_n@>>>E^{w2}_n$ induces a bijection between the set of $G^w_n$-orbits on 
$E^w_n$ and the set of $GL(W^n_{-1})$-orbits on the set of non-degenerate symmetric 
bilinear forms on $W^n_{-1}$. The last set of orbits consists of one element if $p\ne2$
and of two elements (the symplectic forms and the non-symplectic forms) if $p=2$. This 
verifies our claim. 

We see that the first assertion of $\fP_8$ holds.

As above, we identify $E$ with the set of triples $(w,\a,j)$ where $w\in\bE$, $\a$ is a
collection of isomorphisms as in (c) (for each $n\ge0$) and $j$ is a sequence 
$(j_n)_{n\in I}$ where $j_n\in\Bil(W^n_{-1})$ is symmetric non-degenerate.

Assume that $p=2$. Let $J\sub J'\sub I$. From the previous discussion we see that the 
$\bG^\l_0$-orbits on $E$ that contain $E_J$ in their closure and are contained in the
closure of $E_{J'}$ are those of the form $E_K$ where $J\sub K\sub J'$. Let 
$E_{J,J'}=\cup_{K;J\sub K\sub J'}E_K$. We identify $E_J$ with the set of 
$(w,\a,j)\in E$ such that $j_n$ is not symplectic for $n\in J$ and symplectic for 
$n\in I-J$. We identify $E_{J,J'}$ with the set of $(w,\a,j)\in E$ such that $j_n$ is 
not symplectic for $n\in J$ and symplectic for $n\in I-J'$. Let $\tE_J$ be the set of 
all triples $(w,\a,\tj)$ where $w,\a$ are as above and $\tj=(\tj_n)_{n\in I}$ where for
$n\in J$, $\tj_n\in\Bil(W^n_{-1})$ is a symmetric non-symplectic non-degenerate form 
and, for $n\in I-J$, $\tj_n:W^n_{-1}\T W^n_{-1}@>>>\kk$ is the square of a symplectic 
non-degenerate form. 

Now $E_J,E_{J,J'},\tE_J$ are naturally algebraic varieties. Define a finite bijective 
morphism $\s:E_J@>>>\tE_J$ by $(w,\a,j)\m(w,\a,\tj)$ where $\tj_n=j_n$ for $n\in J$, 
$\tj_n=j_n^2$ for $n\in I-J$. Define $\r:E_{J,J'}@>>>\tE_J$ by $(w,\a,j)\m(w,\a,\tj)$ 
where $\tj_n=j_n$ for $n\in J$ and $\tj_n(x,x')=j_n(x,x')^2+j_n(x,x)j_n(x',x')$ for 
$n\in I-J$, $x,x'\in W^n_{-1}$. (To see that this is well defined, we must check that 
for $n\in I-J$, the symplectic form $x,x'\m j_n(x,x')+\sqrt{j_n(x,x)j_n(x',x')}$ on 
$W^n_{-1}$ is non-degenerate. Let $R$ be the radical of this symplectic form. Let 
$H=\{x\in W^n_{-1};j_n(x,x)=0\}$. If $x\in R\cap H$, then $j_n(x,x')$ for all $x'$ 
hence $x=0$. Thus, $R\cap H=0$. Since $H$ is either $W^n_{-1}$ or a hyperplane in 
$W^n_{-1}$, we see that $R\cap H$ is either $R$ or a hyperplane in $R$. It follows that
$\dim R$ is $0$ or $1$. Since $R=\dim W^n_{-1}\mod2$ we see that $\dim R$ is even. 
Hence $R=0$, as required.)

Taking here $J'=I$, we see that $\fP_7$ holds. We now return to a general $J'$. We 
consider the fibre $\cf$ of $\r$ at $(w,\a,\tj)\in\tE_J$. We may identify $\cf$ with 
the set of all collections $(j_n)_{n\in I-J}$ where $j_n\in\Bil(W^n_{-1})$ is symmetric
non-degenerate for all $n$, $j_n$ is symplectic for $n\in I-J'$ and
$\tj_n(x,x')=j_n(x,x')^2+j_n(x,x)j_n(x',x')$ for $n\in I-J$, $x,x'\in W^n_{-1}$. Let 
$\cf'$ be the set of all collections $(h_n)_{n\in I-J}$ where $h_n$ is a linear form 
$W^n_{-1}@>>>\kk$, zero for $n\in I-J'$. We define a map $\cf@>>>\cf'$ by 
$(j_n)_{n\in I-J}\m(h_n)_{n\in I-J}$ where $h_n(x)=\sqrt{j_n(x,x)}$ for 
$x\in W^n_{-1}$. We define a map $\cf'@>>>\cf$ by $(h_n)_{n\in I-J}\m(j_n)_{n\in I-J}$ 
where $j_n(x,x')=\sqrt{\tj_n(x,x')}+h_n(x)h_n(x')$ for $x,x'\in W^n_{-1}$. (We show 
that this is well defined. We must show that $j_n$ given by the last equality is 
non-degenerate. Let $R'$ be the radical of $j_n$. Define $v\in W^n_{-1}$ by 
$h_n(y)=\sqrt{\tj_n(v,y)}$ for all $y\in W^n_{-1}$. If $x\in R',y\in W^n_{-1}$, we have
$\sqrt{\tj_n(x,y)}=h_n(x)h_n(y)=h_n(x)\sqrt{\tj_n(v,y)}$ hence 
$\sqrt{\tj_n(x-h_n(x)v,y)}=0$. Since $\sqrt{\tj_n}$ is non-degenerate we have 
$x-h_n(x)v=0$. Hence $x=h_n(x)v=h_n(h_n(x)v)v=h_n(x)h_n(v)v$. This is $0$ since 
$h_n(v)=\sqrt{\tj_n(v,v)}=0$. Thus $R'=0$.) Clearly, $\cf@>>>\cf',\cf'@>>>\cf$ are 
inverse to each other. We see that $\cf$ is in natural bijection with a vector space of
dimension $\sum_{n\in J'-J}c_n$ where $c_n=\dim W^n_{-1}$. Hence if $\kk,q$ are as in 
$\fP_5$, we have 

$\sum_{K;J\sub K\sub J'}|E_K(\FF_q)|=|E_{J,J'}(\FF_q)|
=\prod_{n\in J'-J}q^{c_n}|E_J(\FF_q)|$.
\nl
From this we see that $|E_K(\FF_q)|=\prod_{n\in K}(q^{c_n}-1)|E_\em(\FF_q)|$ for any 
$K\sub I$. Using this and $\fP_6$ we see that the second assertion of $\fP_8$ holds.

For $\kk,q$ as in $\fP_5$ we have 

$|H^\bla(\FF_q)|=|X^\l(\FF_q)||G(\FF_q)/G^\l_0(\FF_q)|$,
$|X^\l(\FF_q)|=q^{\dim G^\l_3}|E(\FF_q)|$.
\nl
Hence to verify $\fP_5$ it suffices to show that $|E(\FF_q)|$ is a polynomial in $q$ 
with integer coefficients independent of $p$. Using the $\bG^\l_0$-equivariant 
fibration $\ph:E@>>>\bE$ we see that $|E(\FF_q)|=|\bE(\FF_q)||E^w(\FF_q)|$ for any 
$w\in\bE$. Since $|\bE(\FF_q)|$ is a polynomial in $q$ with integer coefficients 
independent of $p$, it suffices to show that for any $w\in\bE(\FF_q)$, $|E^w(\FF_q)|$ 
is a polynomial in $q$ with integer coefficients independent of $p$, or that 
$|E^w_n(\FF_q)|$ is a polynomial in $q$ with integer coefficients independent of $p$ 
for any $w\in\bE(\FF_q)$ and any $n\ge0$. Using the identification 
$E^w_n=E^{w1}_n\T E^{w2}_n$ and the fact that $|E^{w1}_n(\FF_q)|$ is a polynomial in 
$q$ with integer coefficients independent of $p$, we see that it suffices to show that
$|E^{w2}_n(\FF_q)|$ is a polynomial in $q$ with integer coefficients independent of 
$p$. Thus it suffices to check the following statement.

{\it Let $W$ be an $\FF_q$-vector space of dimension $d$. Let $b(W)$ be the set of 
non-degenerate symmetric bilinear forms $W\T W@>>>\FF_q$. Then $|b(W)|$ is a polynomial
in $q$ with integer coefficients independent of $p$.}
\nl
We argue by induction on $d$. For $d=0$ the result is obvious. Assume that $d\ge1$. We 
write $|b(W)|=f(d,q)$. The set of all symmetric bilinear forms $W\T W@>>>\FF_q$ has 
cardinal $q^{d(d+1)/2}$; it is a disjoint union $\sqc_Xb_X(W)$ where $X$ runs over the 
linear subspaces of $W$ and $b_X(W)$ is the set of symmetric bilinear forms 
$W\T W@>>>\FF_q$ with radical equal to $X$. Thus,

$q^{d(d+1)/2}=\sum_X|b_X(W)|=\sum_X|b(W/X)|=\sum_{d'\in[0,d]}g(d,d',q)f(d-d',q)$
\nl
where $g(d,d',q)=|\{X\sub W,\dim X=d'\}|$. We see that

$f(d,q)=q^{d(d+1)/2}-\sum_{d'\in[1,d]}g(d,d',q)f(d-d',q)$.
\nl
Since $g(d,d',q)$ is a polynomial in $q$ with integer coefficients independent of $p$ 
and the same holds for $f(d-d',q)$ with $d'\in[1,d]$ (by the induction hypothesis) it 
follows that $f(d,q)$ is as required.

We see that $\fP_5$ holds.

\head 4. The group $A^1(u)$\endhead
\subhead 4.1\endsubhead
In this section we assume that $p\ge2$ and that $\fP_1$ holds. Let $u\in\cu$. According
to $\fP_1$ there is a unique $\l\in D_G$ such that $u\in X^\l$. Let 
$A^1(u)=Z_{G^\l_1}(u)/Z_{G^\l_1}(u)^0$, a finite $p$-group.

The image of $A^1(u)$ in $Z_G(u)/Z_G(u)^0$ is a normal subgroup (since 
$Z_G(u)=Z_{G^\l_0}(u)$, see 1.1(c), and $Z_{G^\l_1}(u)$ is normal in $Z_{G^\l_0}(u)$).

In this section we describe the finite group $A^1(u)$ in some examples assuming that 
$p=2$ and $G$ is a symplectic group. 

Let $n\ge 1$. Let $I=\{i\in[-n,n];i=n\mod 2\}$. Let $F\in\cc,F\ne0$. Let 
$V=\op_{i\in I}F_i$ where $F_i=F$. Define $N:V@>>>V$ by $(x_i)\m(x'_i)$ where 
$x'_i=x_{i-2}$ for $i\in I-\{-n\},x'_{-n}=0$. We fix $\lar_0\in\Symp(V)$ such that 
$\la(x_i),(y_i)\ra_0=\sum_{i\in I}b(x_i,y_{-i})$ where $b\in\Bil(F)$ satisfies $b^*=b$,
$b$ is non-degenerate and $b\in\Symp(F)$ if $n$ is even.

Let $\lar\in\Symp(V)$ be such that $\la Nx,y\ra+\la x,Ny\ra+\la Nx,Ny\ra=0$ for 
$x,y\in V$ and $\la x,y\ra=\la x,y\ra_0$ if there exists $i$ such that $x_j=0$ for 
$j\ne i$ and $y_j=0$ for $j\ne-i$. We have 
$\la (x_i),(y_i)\ra=\sum_{i,j\in I}b_{ij}(x_i,y_j)$ where $b_{ij}\in\Bil(F)$ are such 
that 

$b_{i-2,j}+b_{i,j-2}+b_{ij}=0$ if $i,j\in I-\{-n\}$,

$b_{i,-i}=b$ for all $i\in I$,

$b_{ii}\in\Symp(F)$ for all $i\in I$, 

$b_{ij}^*=-b_{ji}$ for all $i,j\in I$,
\nl
(We have automatically $b_{ij}=0$ if $i+j\ge 1$.)
\nl
Let $\D'=\{T\in GL(V);TN=NT,\la x,y\ra=\la Tx,Ty\ra\qua\frl x,y\in V\}$, a subgroup of 
$Sp(\lar)$; equivalently, $\D'$ is the set of linear maps $T:V@>>>V$ of the form 

$T:(x_i)\m(x'_i),x'_i=\sum_{j\in I;j\le i}T_{i-j}x_j$
\nl
where $T_r\in\End(F)$ $(r\in\{0,2,4,\do,2n\})$ are such that 
$$b_{ij}(x,y)=\sum_{i',j'\in I;i'\ge i,j'\ge j}b_{i'j'}(T_{i'-i}(x),T_{j'-j}y)
\tag $E_{ij}$ $$
for $i,j\in I,i+j\le0$ and $x,y\in F$. Now $(E_{ij}),(E_{i+2,j-2})$ with $i+j=2k$ are 
equivalent if $(E_{ab})$ is assumed for $a+b=2k+2$ (the sum of those two equations is 
just $E_{i+2,j}$). Thus the conditions that $T$ must satisfy are $E_{ii}$ and 
$E_{i-2,i}$. Setting $x=y$ in these equations we obtain equations $(E_{ii}^0)$,
$(E_{i-2,i}^0)$. Note that the equation $(E_{ii}^0)$ is $0=0$ hence can be omitted; the
equation $(E_{ii})$ is a consequence of $(E_{i-2,i}^0)$ (if it is defined). Hence the 
equations satisfied by the components of $T$ are as follows:
$$(E_{-2,0}^0),(E_{-2,0}),(E_{-4,-2}^0),(E_{-4,-2}),\do,(E_{-n,-n+2}^0),
(E_{-n,-n+2}),(E_{-n,-n})\tag a$$
(for $n$ even),
$$\align&(E_{-1,1}),(E_{-3,-1}^0),(E_{-3,-1}),(E_{-5,-3}^0),(E_{-5,-3}),\do,
(E_{-n,-n+2}^0),\\&(E_{-n,-n+2}),(E_{-n,-n})\tag b\endalign$$ 
(for $n$ odd). Assume first that $n$ is even. The solutions $T_0$ of the first equation
in (a) form an even orthogonal group, a variety with two connected components. For any 
such $T_0$ the solutions $T_2$ of the second equation in (a) form an affine space of 
dimension independent of $T_0$. For any $T_0,T_2$ already determined, the solutions 
$T_4$ of the third equation in (a) form an affine space of dimension independent of 
$T_0,T_2$. Continuing in this way we see that the solutions of the equations (a) form a
variety with two connected components. Moreover, the solutions in which $T_0$ is 
specified to be $1$ form a connected variety.

Assume next that $n$ is odd and $b$ is symplectic. The solutions $T_0$ of the first 
equation in (b) form a symplectic group (a connected variety). For any such $T_0$ the 
solutions $T_2$ of the second equation in (b) form an affine space of dimension 
independent of $T_0$. For any $T_0,T_2$ already determined, the solutions $T_4$ of the 
third equation in (b) form an affine space of dimension independent of $T_0,T_2$. 
Continuing in this way we see that the solutions of the equations (b) form a connected 
variety. Moreover, the solutions in which $T_0$ is specified to be $1$ form a connected
variety.

One can show that, if $n$ is odd, $n\ge3$ and $b$ is not symplectic, the solutions of 
the equations (b) form a variety with two connected components. Moreover, the solutions
in which $T_0$ is specified to be $1$ form a disconnected variety.

In solving the equations above we use repeatedly the statement (c) below. Let $\fQ$ be 
the vector space of quadratic forms $F@>>>\kk$. Define linear maps $a_1,a_2,a_3$ as 
follows: 

$a_1:\End(F)@>>>\fQ(F)$ is $\t\m q,q(x)=b(x,\t(x))$;

$a_2:\{\t\in\End(F);b(\t(x),y)=b(x,\t(y)\qua\frl x,y\in F\}@>>>\Hom(F,\kk)$ is
$\t\m\th,\th(x)=\sqrt{b(x,\t(x))}$; 

$a_3:\{b'\in\Bil(F);b'{}^*=b'\}@>>>\Hom(F,\kk)$ is $b'\m\th,\th(x)=\sqrt{b'(x,x)}$. 
\nl
Then 

(c) $a_1,a_2,a_3$ are surjective.
\nl
For $a_3$ this is clear. Consider now $a_2$. Let $\th\in\Hom(F,\kk)$. By (c) for $a_3$
we can find $b'\in\Bil(F),b'{}^*=b'$ such that $\th(x)=\sqrt{b'(x,x)}$. We can find a 
unique $\t\in End(F)$ such that $b(x,\t(y))=b'(x,y)$. Then $a_2(t)=\th$. Consider now 
$a_1$. Let $q\in\fQ$. Let $b^0$ be a symplectic form on $F$. We can write $b^0=d+d^*$ 
where $d\in\Bil(F)$. We can write $d(x,y)=b(x,\s(y))$ for some $\s\in\End(F)$. Then
$b(x,\s(y))+b(y,\s(x))=b^0(x,y)$. Apply this to the symplectic form 
$b^0(x,y)=q(x+y)+q(x)+q(y)$. Then 

$b(x+y,\s(x+y))+b(x,\s(x))+b(y,\s(y))=b(x,\s(y))+b(y,\s(x))=q(x+y)+q(x)+q(y)$. 
\nl
Hence $b(x,\s(x))+q(x)=\th(x)^2$ for some $\th\in\Hom(F,\kk)$. By (c) for $a_2$ we can 
find $\t\in\End(F)$ such that $b(x,\t(x))=\th(x)^2$. Then $b(x,\s(x))+b(x,\t(x))=q(x)$ 
that is $b(x,(\s+\t)(x))=q(x)$. Thus $a_1$ is surjective. This proves (c). 

\subhead 4.2\endsubhead
Let $V,\lar$ be as in 3.2. Let $N\in\cm_{\lar},V_*=V_*^N$. Let $e$ be as in 2.4. We 
show:

(a) {\it If $W,W'$ are $e$-special subspaces of $V$ (see 2.8) then there exists 
$g\in1+E^{\lar}_{\ge1}V_*$ such that $g(W)=W'$, $gN=Ng$.}
\nl
By 2.8(b) we can find $g_1\in1+E_{\ge 1}V_*$ such that $g_1(W)=W'$, $g_1N=Ng_1$. Then 
$g_1$ carries $\lar$ to a symplectic form $\lar'$ which induces the same symplectic 
form as $\lar$ on $\gr V_*$ and has the same associated quadratic forms as $\lar$ (see 
3.6(b)). By the proof in 3.10 (case 2 and 3) we see that there exists 
$g_2\in1+E_{\ge 1}V_*$ such that $g_2(W')=W'$, $g_2N=Ng_2$ and $g_2$ carries $\lar'$ to
$\lar$. Then $g=g_2g_1$ has the required properties. 

\subhead 4.3\endsubhead
Let $V,\lar,N,V_*,e$ be as in 4.2.

(a) {\it If $\la x,Nx\ra=0$ for any $x\in V_{\ge-1}$, then 
$\cv:=\{g\in E^{\lar}_{\ge1}V_*;gN=Ng\}$ is connected. Hence $A^1(1+N)=\{1\}$.}
\nl
We argue by induction on $e$. Let $\cx$ be the set of all $e$-special subspaces (see 
2.8) of $V$. By 2.8(b) the group $\{g\in 1+E_{\ge 1}V_*;gN=Ng\}$ acts transitively on 
$\cx$. This group is connected (it may be identified as variety with the vector space  
$\{\x\in E_{\ge 1}V_*;\x N=N\x\}$); hence $\cx$ is connected. By 4.2(a), $\cv$ acts 
transitively on $\cx$. Since $\cx$ is connected, it suffices to show that the 
stabilizer $\cv_W$ of some $W\in\cx$ in $\cv$ is connected. This stabilizer is 
$\cv'\T\cv''$ where $\cv',\cv''$ are defined like $\cv$ in terms of $W,W^\pe$ instead 
of $V$. By results in 4.2, $\cv'$ is connected. By the induction hypothesis applied to 
$W^\pe$, $\cv''$ is connected. Hence $\cv'\T\cv''$ is connected. Hence $\cv$ is 
connected. 

\head 5. Study of the varieties $\cb_u$\endhead
\subhead 5.1\endsubhead
We assume that $\kk=\kk_p$. 

We say that an algebraic variety $V$ over $\kk$ has the {\it purity property} if for 
some/any $\FF_q$-rational structure on $V$ (where $\FF_q$ is a finite subfield of 
$\kk$) with Frobenius map $F:V@>>>V$ and any $n\in\ZZ$, any complex absolute value of 
any eigenvalue of $F^*:H^n_c(V,\bbq)@>>>H^n_c(V,\bbq)$ is $q^{n/2}$.

In this section we show that for certain $u\in\cu$ the varieties $\cb_u$ (see 0.1) have
the purity property. We assume that properties $\fP_1-\fP_4,\fP_6,\fP_7$ hold for $G$. 

Let $\l\in D_G$. Let $\Pi^\l$ be the (finite) set of orbits for the conjugation action 
of $G^\l_0$ on $\cb$. Let $\bcb=\{B\in\cb;B\sub G^\l_0\}$. For any $\co\in\Pi^\l$ 
define a morphism $\x^\co:\co@>>>\bcb$ by $B\m(B\cap G^\l_0)G^\l_1$. We show:

(a) {\it The fibres of $\x^\co:\co@>>>\bcb$ are exactly the orbits of $G^\l_1$ acting 
on $\co$ by conjugation.}
\nl
If $B,B'\in\co$, $\x^\co(B)=\x^\co(B')$, then $B'=g\i Bg$ with $g\in G^\l_0$, 
$(B'\cap G^\l_0)G^\l_1=(B\cap G^\l_0)G^\l_1=g\i(B\cap G^\l_0)G^\l_1g$. Hence 
$g\in(B\cap G^\l_0)G^\l_1$. Writing $g=g'g'',g'\in B\cap G^\l_0$, $g''\in G^\l_1$, we 
have $B'=g\i Bg=g''{}\i Bg''$. This proves (a).

Let $Y^\l=\{(u,B)\in X^\l\T\cb;u\in B\}$. We have a partition 
$Y^\l=\cup_{\co\in\Pi^\l}Y^\l_\co$ where $Y^\l_\co=\{(u,B)\in X^\l\T\co;u\in B\}$. Let 
$\co\in\Pi^\l$. We show:

(b) {\it$Y^\l_\co$ is smooth.}
\nl
Let $\tB\in\co$. Let $Y'=\{(u,g)\in X^\l\T G^\l_0;g\i ug\in\tB\cap X^\l\}$. We have a
fibration $Y'@>>>Y^\l_\co$ with smooth fibres isomorphic to $G^\l_0\cap\tB$. Hence it 
suffices to show that $Y'$ is smooth. Let $Y''=(\tB\cap X^\l)\T G^\l_0$. Define
$Y'@>\si>>Y''$ by $(u,g)\m(g\i ug,g)$. It suffices to show that $Y''$ is smooth, or 
that $\tB\cap X^\l$ is smooth. But $\tB\cap X^\l$ is open in $\tB\cap G^\l_2$ which is 
smooth, being an algebraic group. This proves (b).

For any $\b\in\bcb$ let $\cg_\b^\co=((B\cap G^\l_2)G^\l_3)/G^\l_3$ where $B\in\co$ is 
such that $\x^\co(B)=\b$. Note that $\cg_\b^\co$ is a closed connected subgroup of 
$G^\l_2/G^\l_3$, independent of the choice of $B$. (To verify the last statement it  
suffices, by (a), to show that for $B$ as above and $v\in G^\l_1$ we have 
$(vBv\i\cap G^\l_2)G^\l_3=(B\cap G^\l_2)G^\l_3$. This follows from 1.1(b).) Now 
$G^\l_0$ acts on $\bcb$ and on $G^\l_2/G^\l_3$ by conjugation. From the definitions we 
see that for $g\in G^\l_0$ and $\b\in\bcb$ we have $\cg_{g\b g\i}^\co=g\cg_\b^\co g\i$.

Let $\bY^\l_\co=\{(x,\b)\in\bX^\l\T\bcb;x\in\cg_\b^\co\}$. We show:

(c) {\it$\bY^\l_\co$ is a closed smooth subvariety of $\bX^\l\T\bcb$.}
\nl
Let $\tB\in\co$. We have a fibration $X^\l\T G^\l_0@>>>\bX^\l\T\bcb$, 
$(u,g)\m(\p^\l(u),\x^\co(g\tB g\i))$ with smooth fibres. It suffices to show that the
inverse image of $\bY^\l_\co$ under this fibration is a closed smooth subvariety of 
$X^\l\T G^\l_0$, or that 

$\{(u,g)\in X^\l\T G^\l_0;g\i ug\in X^\l\cap((\tB\cap G^\l_2)G^\l_3)\}$
\nl
is a closed smooth subvariety of $X^\l\T G^\l_0$, or that 
$(X^\l\cap((\tB\cap G^\l_2)G^\l_3))\T G^\l_0$ is a closed smooth subvariety of 
$X^\l\T G^\l_0$ or that $X^\l\cap((\tB\cap G^\l_2)G^\l_3)$ is a closed smooth 
subvariety of $X^\l$. It is closed since $(\tB\cap G^\l_2)G^\l_3$ is closed in 
$G^\l_2$. It is smooth since it is an open subset of $(\tB\cap G^\l_2)G^\l_3$ which is 
smooth, being an algebraic group.

We show:

(d) {\it The morphism $a:Y^\l_\co@>>>\bY^\l_\co,(u,B)\m(\p^\l(u),\x^\co(B))$ is a 
fibration with fibres isomorphic to an affine space of a fixed dimension.}
\nl
Clearly, $a$ is surjective. Let $(u,B)\in Y^\l_\co$. Let

$Z:=a\i(a(u,B))
=\{(u',B');u=u'f,B'=vBv\i\text{ for some }v\in G^\l_1,f\in G^\l_3;u'\in B'\}$.
\nl
We show only that $Z$ is isomorphic to an affine space of fixed dimension. Let 
$\tZ=\{(f,v)\in G^\l_3\T G^\l_1;v\i uf\i v\in B\}$. Then $Z=\tZ/(B\cap G^\l_1)$ where 
$B\cap G^\l_1$ acts freely on $\tZ$ by $b:(f,v)\m(f,vb\i)$. Since conjugation by
$G^\l_1$ acts trivially on $G^\l_2/G^\l_3$, the map $(f,v)\m(f',v),f'=u\i v\i uf\i v$ 
is an isomorphism 
$$\align&\tZ@>>>\tZ'=\{(f',v)\in G^\l_3\T G^\l_1;uf'\in B\}=
\{(f',v)\in G^\l_3\T G^\l_1;f'\in B\}\\&=(G^\l_3\cap B)\T G^\l_1\endalign$$
(we use $u\in B$) and we have $Z=(G^\l_3\cap B)\T G^\l_1/(B\cap G^\l_1)$. Now 
$G^\l_3\cap B,G^\l_1,B\cap G^\l_1$ are connected unipotent groups of dimension
independent of $B$, for $B\in\co$. (The connectedness follows from the fact that these
unipotent groups are normalized by a maximal torus of $G$ contained in $G^\l_0\cap B$. 
The fact that the dimension does not depend on $B$ follows from the fact that 
$G^\l_1,G^\l_3$ are normalized by $G^\l_0$.) We see that $Z$ is an affine space of 
constant dimension.

We now fix $x\in\bX^\l$. Let $\Si=(\p^\l)\i(x)\sub X^\l$. Let 
$\cb_\Si=\{(u,B)\in\Si\T\cb;u\in B\}$. We have $\cb_\Si=\sqc_{\co\in\Pi^\l}\co_\Si$ 
where $\co_\Si=\{(u,B)\in\Si\T\co;u\in B\}$. Let $\co\in\Pi^\l$. Let 
$\bcb^\co_x=\{\b\in\bcb;x\in\cg^\co_\b\}$. We show:

(e) {\it$\bcb^\co_x$ is a closed subvariety of $\bcb$ and $a':\co_\Si@>>>\bcb^\co_x$, 
$(u,B)\m\x^\co(B)$ is a fibration with fibres isomorphic to an affine space of a fixed 
dimension.}
\nl
Let $\tB\in\co,u_0\in\Si$. We have a locally trivial fibration $G^\l_0@>>>\bcb$, 
$g\m\x^\co(g\tB g\i)$. To show that $\bcb^\co_x$ is closed it suffices to show that its
inverse image under this fibration is closed in $G^\l_0$, or that 
$\{g\in G^\l_0;g\i u_0g\in(\tB\cap G^\l_2)G^\l_3\}$ is closed in $G^\l_0$. This is 
clear since $(\tB\cap G^\l_2)G^\l_3$ is closed in $G^\l_2$. The second assertion of (e)
follows from (d) using the cartesian diagram
$$\CD
\co_\Si@>a'>>\bcb^\co_x\\
@VVV         @VVV\\
Y^\l_\co@>a>>\bY^\l_\co
\endCD$$
where the left vertical map is the obvious inclusion and the right vertical map is
$\b\m(x,\b)$.

(f) {\it If the closure of the $G^\l_0$-orbit in $G^\l_2$ of some/any $u\in\Si$ is a 
subgroup $\G$ of $G^\l_2$ then $\bcb^\co_x$ is smooth.}
\nl
Let $\tB\in\co,u_0\in\Si$. As in the proof of (e) it suffices to show that
$\{g\in G^\l_0;g\i u_0g\in(\tB\cap G^\l_2)G^\l_3\}$ is smooth. This variety is a
fibration over $R=(G^\l_0-\text{conjugacy class of }u_0)\cap((\tB\cap G^\l_2)G^\l_3)$
with smooth fibres isomorphic to $Z_{G^\l_0}(u_0)$ (via $g\m g\i u_0g$). Hence it 
suffices to show that $R$ is smooth. From our assumption we see that $R$ is open in
$\G\cap((\tB\cap G^\l_2)G^\l_3)$ which is smooth being an algebraic group. This proves
(f).

\mpb

Note that the hypothesis of (f) holds at least in the case where the $G^\l_0$-conjugacy
class of some/any $u\in\Si$ is open dense in $G^\l_2$. We show

(g) {\it If the hypothesis of (f) holds then $\cb_\Si$ has the purity property.}
\nl
From (e),(f) we see that $\bcb^\co_x$ is a smooth projective variety of pure dimension.
From \cite{\DE} it then follows that $\bcb^\co_x$ has the purity property. From this 
and (e) we see that for $\co\in\Pi^\l$, $\co_\Si$ has the purity property. Using this 
and the partition $\cb_\Si=\sqc_{\co\in\Pi^\l}\co_\Si$, we see that (g) holds.

\subhead 5.2\endsubhead
Let $\bZ(x)=\{\bg\in\bG^\l_0;\bg x=x\bg\}$. Let $\tZ(x)$ be the inverse image of 
$\bZ(x)$ under $G^\l_0@>>>\bG^\l_0$. Thus we have $G^\l_1\sub\tZ(x)$ and 
$\tZ(x)/G^\l_1@>\si>>\bZ(x)$. Note that the inverse image of $\bZ(x)^0$ is $\tZ(x)^0$ 
and we have $\tZ(x)^0/G^\l_1@>\si>>\bZ(x)^0$. Now $\tZ(x)$ acts transitively (by 
conjugation) on $\Si$. (Indeed, let $u,u'\in\Si$. By $\fP_6$ we can find $g\in G^\l_0$
such that $u'=gug\i$. We have automatically $g\in\tZ(x)$.) Since $\Si$ is irreducible, 
it follows that $\tZ(x)^0$ acts transitively (by conjugation) on $\Si$. 

Let $u\in\Si$. Recall that $\cb_u=\{B\in\cb;u\in B\}$. Let 
$Z'_G(u)=Z_G(u)\cap\tZ(x)^0$. Since $Z_G(u)\sub\tZ(x)$, see 1.1(c), we see that 
$Z'_G(u)$ is a normal subgroup of $Z_G(u)$ containing $Z_G(u)^0$. Let $A'(u)$ be the 
image of $Z'_G(u)$ in $A(u):=Z_G(u)/Z_G(u)^0$; this is a normal subgroup of $A(u)$. We 
have $Z'_G(u)/Z_{G^\l_1}(u)@>\si>>\bZ(x)^0$. Hence
$Z'_G(u)=Z_{G^\l_1}(u)Z'_G(u)^0=Z_{G^\l_1}(u)Z_G(u)^0$. It follows that

{\it$A'(u)$ is the image of the obvious homomorphism} $A^1(u)@>>>A(u)$. 
\nl
Now $Z_G(u)$ acts by conjugation on $\cb_u$; this induces an action of $A(u)$ on 
$H^n(\cb_u,\bbq)$ which restricts to an $A'(u)$-action on $H^n(\cb_u,\bbq)$.

Assume that the hypothesis of 5.1(f) holds and that $A'(u)$ acts trivially 
on $H^n_c(\cb_u,\bbq)$ for any $n$. We show:

(a) {\it$\cb_u$ has the purity property.}
\nl
Define $f:\cb_\Si@>>>\Si$ by $(g,B)\m g$. For any $n$, $R^nf_!(\bbq)$ is an equivariant
constructible sheaf for the transitive $\tZ(x)^0$ action on $\Si$; hence it is a local 
system on $\Si$ corresponding to a representation of $A'(u)$ (the group of components 
of the isotropy group of $u$ in $\tZ(x)^0$) on $H^n_c(\cb_u,\bbq)$. This representation
is trivial hence $R^nf_!(\bbq)$ is a constant local system. Since $\Si$ is an affine 
space of dimension say $d$ we see that $H^a_c(\Si,R^nf_!(\bbq))$ is 
$H^n_c(\cb_u,\bbq)(-d)$ if $a=2d$ and is zero if $a\ne2d$. It follows that the standard
spectral sequence 

$E_2^{a,n}=H^a_c(\Si,R^nf_!(\bbq))\Rightarrow H^{a+n}_c(\cb_\Si,\bbq)$
\nl
is degenerate. Hence the purity property of $\cb_\Si$ (see 5.1(g)) implies that any 
complex absolute value of any eigenvalue of the Frobenius map on

$E_2^{2d,n}=H^n_c(\cb_u,\bbq)(-d)$
\nl
is $q^{d+n/2}$. Hence any complex absolute value of any eigenvalue of the Frobenius map
on $H^n_c(\cb_u,\bbq)$ is $q^{n/2}$. This proves (a).

\subhead 5.3\endsubhead
Since the hypothesis of 5.1(f) is not satisfied in general, we seek an alternative way 
to prove purity.

Let $\g$ be the $\bG^\l_0$-orbit of $x$ in $\bX^\l$. Let $\hag@>\r>>\g_1@<\s<<\g$ be as
in $\fP_7$. Let $\Xi=\r\i(\s(x))$. Let 
$\bcb^\co_\Xi=\{(x',\b)\in\bY^\l_\co;x'\in\Xi\}$, a closed subvariety of $\bY^\l_\co$. 
We show:

(a) {\it$\bcb^\co_\Xi$ is smooth of pure dimension.}
\nl
Let $\b_0\in\bcb$. Let $\cg_0=\cg^\co_{\b_0}$. It suffices to show that the inverse 
image of $\bcb^\co_\Xi$ under the fibration $\Xi\T\bG^\l_0@>>>\Xi\T\bcb$,
$(x',\bg)\m(x',\bg\b_0\bg\i)$ (with smooth connected fibres) is smooth of pure 
dimension, or that $\fS:=\{(x',\bg)\in\Xi\T\bG^\l_0;\bg\i x'\bg\in\cg_0\}$ is smooth of
pure dimension. The morphism $f:\fS@>>>\hag\cap\cg_0$, $(x',\bg)\m\bg\i x'\bg$ is 
smooth with fibres of pure dimension. (We show only that for any $y\in\hag\cap\cg_0$,
the fibre $f\i(y)$ is isomorphic to $\{\bg\in\bG^\l_0;\bg x\bg\i=x\}$ which is smooth
of pure dimension. We have
$$\align&f\i(y)=\{(x',\bg)\in\Xi\T\bG^\l_0;\bg\i x'\bg=y\}\cong
\{\bg\in\bG^\l_0;\bg y\bg\i\in\Xi\}\\&=\{\bg\in\bG^\l_0;\r(\bg y\bg\i)=\s(x)\}=
\{\bg\in\bG^\l_0;\bg\s\i(\r(y))\bg\i=x\}\endalign$$
and it remains to use the transitivity of the $\bG^\l_0$-action on $\g$.) It suffices 
to show that $\hag\cap\cg_0$ is empty or smooth, connected. Now $\hag$ is open in 
$G^\l_2/G^\l_3$ hence $\hag\cap\cg_0$ is open in $\cg_0$ which is connected and smooth 
(being an algebraic group).

We show:

(b) {\it Assume that for any $\co\in\Pi^\l$ there is a $\kk^*$-action on $\bcb^\co_\Xi$
which is a contraction to the projective subvariety $\bcb^\co_x$. Then $\cb_\Si$ has 
the purity property.}
\nl
Consider an $\FF_q$-rational structure on $G$ such that $G^\l_a$ is defined over 
$\FF_q$ for any $a$ and $\co,x,\Xi$ are defined over $\FF_q$. Let $\z$ be an eigenvalue
of Frobenius on $H^n(\bcb^\co_x,\bbq)$. By \cite{\DEII, 3.3.1}, any complex absolute 
value of $\z$ is $\le q^{n/2}$ (since $\bcb^\co_x$ is projective). Our assumption 
implies that the inclusion $\bcb^\co_x\sub\bcb^\co_\Xi$ induces for any $n$ an 
isomorphism $H^n(\bcb^\co_\Xi,\bbq)@>\si>>H^n(\bcb^\co_x,\bbq)$. Hence $\z$ is also 
an eigenvalue of Frobenius on $H^n(\bcb^\co_\Xi,\bbq)$. Since $\bcb^\co_\Xi$ is smooth 
of pure dimension say $d$, it satisfies Poincar\'e duality; hence $q^d\z\i$ is an 
eigenvalue of Frobenius on $H^{2d-n}_c(\bcb^\co_\Xi,\bbq)$. By \cite{\DEII, 3.3.1} 
applied to $\bcb^\co_\Xi$, we see that any complex absolute value of $q^d\z\i$ is 
$\le q^{(2d-n)/2}$ hence any complex absolute value of $\z$ is $\ge q^{n/2}$. It 
follows that any complex absolute value of $\z$ is $q^{n/2}$. We see that $\bcb^\co_x$ 
has the purity property. (This argument is similar to one of Springer in \cite{\SP}.) 
From this and 5.1(e) we see that for $\co\in\Pi^\l$, $\co_\Si$ has the purity property.
Using this and the partition $\cb_\Si=\sqc_{\co\in\Pi^\l}\co_\Si$, we see that 
$\cb_\Si$ has the purity property.

\mpb

If we assume in addition that $A'(u)$ acts trivially on $H^n_c(\cb_u,\bbq)$ for any $n$
we see as in 5.2 that $\cb_u$ has the purity property.

\subhead 5.4\endsubhead
Let $V,\lar$ be as in 3.2. Assume that $p=2$ and that $G=Sp(\lar)$. Let $u\in\cu$. We 
set $u=1+N,V_*=V_*^N$. Assume that 

(a) $\la x,Nx\ra=0$ for any $x\in V_{\ge-1}$. 
\nl
We set 
$$\G=1+\{N'\in E_{\ge2}^{\lar}V_*;\la x,N'x\ra=0\qua\frl x\in V_{\ge-1}\}.$$
Now $\G$ is a subgroup of $1+E_{\ge2}^{\lar}V_*$. (Assume that $1+N',1+N''\in\tG_2$.
Let $x\in V_{\ge-1}\}$. We have $\la x,N'x\ra=0,\la x,N''x\ra=0$. We must show that 
$\la x,(N'+N''+N'N'')x\ra=0$ or that $\la x,N'N''x\ra=0$. This follows from 
$N'N''x\in V_{\ge3}$ and $3-1\ge1$.) Clearly, $\G$ is normal in $G^\l_0$. Since $\G$ is
a closed unipotent subgroup normalized by $G^\l_0$, it must be connected. Now 

$\cj:=1+\{N'\in\tG_2;\bN'\in\End_2^0(\gr V_*)\}$
\nl
is open in $\G$ since it is the inverse image under $\G@>>>\End_2(\gr V_*),1+N'\m\bN'$ 
of the open subset $\End_2^0(\gr V_*)$ of $\End_2(\gr V_*)$. Also $\cj\ne\em$ since 
$1+N\in\cj$. Hence $\cj$ is an open dense subset of $\G$. By results in 3.14, $\cj$ is 
the $G^\l_0$-conjugacy class of $1+N$. 

We see that the hypothesis of 5.1(f) holds. Using 5.2(a) we see that:

(b) {\it$\cb_u$ has the purity property for any $u\in G$ whose conjugacy class is 
minimal in the unipotent piece containing it, see 1.1, and such that any Jordan block 
of even size appears an even number times.}
\nl
(For such $u$, $A'(u)$ is trivial by 4.3(a).)

Alternatively, one can show that for $u$ as in (b) the method of 5.3 is applicable (the
hypothesis of 5.3(b) holds) and one obtains another proof of (b).

\enddocument